\documentclass[review]{elsarticle}

\usepackage{lineno,hyperref}
\modulolinenumbers[5]
\usepackage{enumitem}
\usepackage{array}
\usepackage{tabularx} 
\usepackage{adjustbox} 
\usepackage{makecell}
\usepackage{amsmath}
\usepackage{float}
\usepackage{graphicx}
\usepackage{subfigure}
\usepackage{geometry}
\journal{Journal of \LaTeX\ Templates}

\begin{document}
\nolinenumbers

\begin{frontmatter}

\title{Perfectionism Search Algorithm (PSA): An Efficient Meta-Heuristic
Optimization Approach}

\author[1]{Amin Ghodousian\corref{mycorrespondingauthor}}
\cortext[mycorrespondingauthor]{Corresponding author}
\ead{a.ghodousian@ut.ac.ir}
\author[2]{Mahdi Mollakazemiha}
\author[3]{Noushin Karimian}

\address[1]{School of Engineering Science, College of Engineering, University of Tehran, Tehran, Iran.}

\address[2]{Faculty of Mathematical Sciences, Shahid Beheshti University, Tehran, Iran.}

\address[3]{Department of Electrical and Electronic Engineering, The University of Manchester, Manchester, UK}

\begin{abstract}

This paper proposes a novel population-based meta-heuristic optimization algorithm, called Perfectionism Search
Algorithm (PSA), which is based on the psychological aspects of perfectionism. The PSA
algorithm takes inspiration from one of the most popular model of perfectionism, which was proposed by Hewitt and
Flett. During each iteration of the PSA algorithm, new solutions are generated by
mimicking different types and aspects of perfectionistic behavior. In order to have a complete perspective on the
performance of PSA, the proposed algorithm is tested with various nonlinear optimization problems, through
selection of 35 benchmark functions from the literature. The generated solutions for these
problems, were also compared with 11 well-known meta-heuristics which had been applied to many complex and
practical engineering optimization problems. The obtained results confirm the high performance of the proposed
algorithm in comparison to the other well-known algorithms. 

\end{abstract}


\end{frontmatter}

\linenumbers

\section{Introduction}
\nolinenumbers

\noindent
Over the past decades, development and application of optimization models have attracted growing attention
amongst mathematicians and engineers. In recent years, , they have received enormous attention, due to the rapid
progress of computer technology, including development and availability of different types of user-friendly software,
high-speed and parallel processors, and artificial neural networks.
Even though most practical optimization problems have some restrictions that must be satisfied, study of available
techniques for unconstrained optimization problems is important for several reasons. For example, many algorithms
solve a constrained problem by converting it into a sequence of unconstrained problems via either Lagrangian
multipliers or penalty and barrier functions. Moreover, most methods proceed by finding an improving direction and
then function minimization takes place along this direction. This line search is equivalent to minimizing a function of
one variable without constraints or with simple constraints, such as lower and upper bounds of the variables [2]. It is worth stating that, improving directions are the backbone of many well-known optimization algorithms such as the
Cyclic Coordinate method, the Hooke and Jeeves, the Rosenbrock, the Steepest Descent method, the Newton’s
method and its modifications, and Quasi-Newton and Conjugate Gradient methods.\\

Although, the above-mentioned mathematical based methods are fast, quick in finding a local minimum, and can
guarantee global optima in both simple and ideal models, but they suffer from limitations, as the optimality
conditions and convergence of these procedures all rely on some specific assumptions. For instance, some types of
convexity, the first (or higher) order differentiability of an objective function; continuity, linearity or the order of
nonlinearity of constraints; the linear independence or orthogonality of directions along which an objective function
is optimized and the closeness of some mappings or regions in problems. Furthermore, to guarantee convergence to
the global optimum, some algorithms require careful attention in selecting the initial values. However, these
assumptions do not hold true for many problems, and in any case, they cannot be easily verified. \\
On the other hand, over the past decades, there has been a growing interest in the development of meta-heuristic
algorithms. These algorithms are inspired by the behaviors of natural phenomena and do not require any
consideration of previous prerequisites [31]. Thus, meta-heuristics can allow for near-optimum solutions to be found
within a reasonable computation-time, with efficient use of memory, without any loss of subtle nonlinear
characteristics of the model and any need for complex derivatives or gradient information, and require careful
selection of initial values. Metaheuristic optimization algorithms rely on rather simple concepts and are easy to
implement. Moreover, these algorithms can bypass local optima and applied to a wide range of problems covering
different disciplines. Additionally, it has been experimentally shown that these algorithms which are able to perform
well, can provide suitable solutions for complex engineering optimization problems. These algorithms are in fact,
good alternative approach to traditional methods, especially for solving combinatorial problems, NP-hard problems
and problems in which the search space grows exponentially with the size of the problem [18]. For all these reasons,
many meta-heuristics were designed to solve various types of applied problems such as Scheduling problems [56],
pattern recognition [15], data clustering [37], tuning of neural networks [33], data mining [64], engineering [5,20]
and optimization [4,22,36,54,55]. Also, in recent years, there has been a growing attempt in developing algorithms
inspired by nature, providing some proofs of convergence, solving different classes of optimization problems and
investigating some other theoretical analyses [11,18]. \\
Nature-inspired metaheuristic algorithms can be roughly classified into four categories according to their
inspirations: evolution-based, physics-based, swarm-based and human-based methods. Evolution-based methods are
inspired by the laws of natural evolution. The most popular evolution-inspired technique is Genetic Algorithm [27].
Other popular algorithms are Biogeography Based Optimization (BBO) [44], Evolution Strategy (ES) [42] and
Probability Based Incremental Learning (PBIL) [9]. Physics-based methods imitate the physical rules in the universe.
Examples of these algorithms include Gravitational Search Algorithm (GSA) [41], Black Hole (BH) algorithm [22],
Simulated Annealing [32], Central Force Optimization (CFO) [14] and Water wave optimization (WWO) [65]. The
third group of nature-inspired methods includes swarm-based techniques that imitate the social behavior of groups of
animals. The most popular algorithms are Particle Swarm Optimization (PSO) [29] and Ant Colony Optimization
(ACO) [10,48]. Other examples of swarm-based techniques are Bees Algorithm (BA) [38], Firefly Algorithm (FA)
[61], Whale Optimization Algorithm (WOA) [35], Forest Optimization Algorithm (FOA) [17], Social Spider
Optimization (SSO) [8], Cuckoo Optimization Algorithm (COA) [39], Grey Wolf Optimizer (GWO) [34], Lion Pride
Optimizer (LPO) [58] and Bat-inspired Algorithm (BA) [60]. The final category consists of methods which are
inspired by human behaviors. Some of the most popular algorithms are Teaching Learning Based Optimization
(TLBO) [40], Harmony Search (HS) [16], Cultural Algorithm (CA) [43], Group Search Optimizer (GSO) [23],
Exchange Market Algorithm (EMA) [19] and Group Counseling Optimization (GCO) algorithm [12]. \\
In addition to the fundamental differences between these algorithms in terms of their inspirations, they also differ in
many other ways such as: types of selection strategies, the generation of new solutions and the definition of
exploitation and exploration operators. For instance, CFO, Stochastic Algorithm (SA) and (1+1)- Evolution Strategy
(ES) [30] are the three approaches that are shown to be very different amongst the heuristic algorithms. CFO is a
deterministic heuristic search algorithm which works based on the metaphor of gravitational kinematics and does not
use any random parameter in its formulation; SA is a stochastic algorithm in which the search starts from a single point and continues in a sequential manner; and (1+1) - ES is a simple evolutionary algorithm with one parent
generating one offspring per iteration. However, most of the heuristic algorithms have a stochastic behavior, where
they do search in a parallel manner with multiple initial points and generate several offsprings from several parents.
In spite of the previously mentioned differences among the heuristic algorithms, population-based meta-heuristic
optimization algorithms share a common feature regardless of their nature. The search process is divided into two
phases: (i) exploration - the ability to expand the search space, and (ii) exploitation - the ability to find the optima
around a good solution [1]. Finding a right balance between exploration and exploitation is the most challenging task
in the development of any metaheuristic algorithm. Moreover, the individuals within population-based algorithms
usually pass three steps to complete the two phases of exploration and exploitation: (i) self-adaptation (each
individual improves its performance), (ii) cooperation (individuals cooperate with each other through transfer of
information), and (iii) competition (individuals compete to survive).\\
It is important to mention that there is no single algorithm to achieve the most appropriate solutions for all classes of
problems. Some algorithms, in comparison to others, provide better solutions for some particular problems.
Generally, the effectiveness of different optimization algorithms vary based on their generality, reliability, precision,
sensitivity to parameters and data, computational effort and convergence. Therefore, pursuing for a new optimization
technique based on its similarities with either natural or artificial phenomena is an open problem [59].
This study describes a novel human-based metaheuristic optimization algorithm (namely, Perfectionism Search
Algorithm - PSA). The source of inspiration for PSA is the behaviors of extreme perfectionists from a psychological
perspective. Perfectionism is a multi-dimensional personality trait characterized by a person’s continuous striving for
flawlessness and setting exceedingly high standards of performance, accompanied by tendencies for overly critical
evaluations of one’s own behavior [52]. Perfectionism is best conceptualized as a multi-dimensional characteristic
[51], since psychologists agree that such type of characteristic consists of many positive and negative aspects. The
maladaptive form of perfectionism, drives people to try and achieve an unattainable ideal, while the adaptive form of
perfectionism can motivate them to reach their goals. When perfectionists do not reach their goals, they often fall
into depression [62]. Recognizing that perfectionism has personal and social dimensions, Hewitt and Flett [24]
proposed a model with an ability to distinguish between the three forms of perfectionism: Self-Oriented
perfectionism (SOP), Other-Oriented perfectionism (OOP), and Socially Prescribed Perfectionism (SPP), where the
three forms comprise of different attitudes, motivations, and behaviors. \\
The model proposed by Hewitt and Flett which explains the concept of perfectionism has been the main source of
inspiration for the development of PSA algorithm. In the proposed algorithm, the above-mentioned types of
perfectionism are mathematically modeled. Then, PSA assigns a probability to each type, which are dynamically
changed by each round of iterations. During each iteration, PSA attempts to generate perfect solutions; i.e. solutions
with higher objective function values. In order to generate a new solution, the algorithm initially selects one type of
perfectionism based on its associated probability, and then generates a solution according to the behaviors of
perfectionists belonging to that particular group. In each case, a comparison is made to decide if the new generated
solution is a good solution or not (compared with the worst solution in the population). If the generated solution is
not a good one, the probability of the selected type of perfectionism is decreased for the next round of generation.
This case is correlating to the depression case when perfectionists do not reach their goals. The efficiency of the PSA
algorithm proposed and developed in this research is evaluated by solving 35 mathematical optimization problems.
Optimization results from these problems demonstrate that PSA is an efficient algorithm for complex optimization
problems compared to the 11 well-known metaheuristic algorithms. \\
The remainder of this paper is organized as follows. In section 2, some characteristics of SOP, OOP and SPP
perfectionists are briefly discussed. This section provides some background information which are the psychological
sources of inspiration for the proposed algorithm. In Section 3, the proposed algorithm (PSA) is outlined, and its
characteristics and implementation steps are explained in details. Section 4 provides a discussion of the proposed
approach which presents the most considered theoretical aspects and applied suggestions concerning PSA.
Comparative study, test problems and experimental results are presented and discussed in Section 5. In this section,
the performance of the presented algorithm is illustrated with several benchmark datasets and compared with ACO,
PSO, TLBO, BBO, Differential Evolution (DE) [53], HS, FA, WOA, BA, GSA and BH. Finally, Section 6 summarizes the main findings of this study and suggests directions for future research.

\section{Basic Rules for PSA – A Brief Conceptual Discussion on Perfectionism }
As mentioned in the previous section, Hewitt and Flett’s tripartite perfectionism model divides perfectionism into
three forms of SOP, OOP and SPP. In the following, a brief overview of important characteristics for these three
forms will be provided. This discussion provides a complete background into creating an outline of PSA, which is
done by deriving the basic aspects from the perfectionistic behaviors. In the next section, these basic rules are
mathematically modeled in order to perform the optimization.

\subsection{Self-Oriented Perfectionism (SOP)}
SOP reflects beliefs that are striving for perfection and generally, being perfect is important. Self-oriented
perfectionists have exceedingly high personal standards from themselves, continuous expectation of being perfect,
and are very self-critical if they fail to meet the expected criteria and demands [24]. In addition to the above, SOP
can be considered to be double-edged, in a way that can also be associated to some positive characteristics such as
conscientiousness, nurturing, intimacy, social development, and altruism [49,51]. On the other hand, these
characteristics often cause self-esteem deficits and self-evaluation and make individuals more prone to depression
[46].

\subsection{Other-Oriented Perfectionism (OOP)}
OOP reflects beliefs, where it is important for others to strive for perfection and be perfect [24,26]. Other-oriented
perfectionists have exceedingly high standards from others, expect others to be perfect, and are also highly critical of
others who fail to meet the expectations [24,26]. So, individuals with tendencies towards OOP have behaviors that
are similar to self-oriented perfectionists, with a difference of behavior directed outward towards others rather than
requiring the self to be perfect. Such characteristics result in lack of trust towards others [24]. OOP is similarly
double-edged and displays positive associations with grandiose narcissism, Machiavellianism, psychopathy, and
aggressive humor [45, 49,50]. Nonetheless, OOP is also tied to lower burnout [7], superior problem solving [13], and
positive self-regard [49].

\subsection{Socially Prescribed Perfectionism (SPP)}
Socially prescribed perfectionists believe that exceedingly high standards are being imposed on them. They
believe others expect them to be perfect, and think that others will be highly critical of them if they fail to meet their
expectations [24,26]. SPP involves a tendency towards having a fear of negative social evaluation and desire for
approval from those around them. When the perceived standards are not met, individuals blame themselves and
lower their self-worth, consequently resulting in depressive symptoms [24]. Socially prescribed perfectionists think
and behave in ways that generate stress, which in turn increases the risk for depression disorder [25].\\
For ease of future reference and before going into in-depth description of the PSA algorithm, the basic rules of
PSA algorithm, which have been derived from the above statements, will be summarized below. The generated
solution from PSA algorithm applied to an optimization problem, is based on the following two main rules:\\
Rule1 (Striving step) – Perfectionists strive for flawlessness, set high standards of performance and have
tendencies for overly critical evaluations of one’s behavior [52].\\
Rule2 (Depression step) – When perfectionists do not reach their expectation, they often fall into depression mood
[62]. \\
In PSA, implementation of the striving step is conducted by modeling some basic characteristics of the above
three types of perfectionism. These characteristics can be expressed as follows:\\

\scriptsize CHARACTERISTICS OF SELF-ORIENTED PERFECTIONISTS IN PSA\\
\begin{tabular}{|@{}|l|l} 
 \hline
 \textbf{SOP1} & They expect to be perfect [24]. \\ 
 \hline
 \textbf{SOP2} & They have exceedingly high personal standards [24]. \\ 
 \hline
\end{tabular}
\\
\\
\scriptsize CHARACTERISTICS OF OTHER-ORIENTED PERFECTIONISTS IN PSA\\
\begin{tabular}{|@{}|l|l} 
 \hline
 \textbf{OOP1} & They expect others to be perfect [24,26]. \\ 
 \hline
 \textbf{OOP2} & They have high standards from others [24,26]. \\ 
 \hline
\textbf{OOP3} & They have lack of trust towards others [24,26]. \\ 
 \hline
 \textbf{OOP4} & They display a positive correlation with narcissism [45,49,50]. \\ 
 \hline
\end{tabular}
\\
\\
\scriptsize CHARACTERISTICS OF SOCIALLY PRESCRIBED PERFECTIONISTS IN PSA\\
\begin{tabular}{|@{}|l|l} 
 \hline
 \textbf{SPP1} & They believe others expect them to be perfect [24]. \\ 
 \hline
 \textbf{SPP2} & They desire approval from those around them [24]. \\ 
 \hline
\end{tabular}
\\
\\
Also, the depression step is performed by modeling the rate of depression based on the following observations:\\
\scriptsize DEPRESSION EFFECTS\\
\begin{tabular}{|@{}|l|l} 
 \hline
 \textbf{Dep1} & \makecell{When perceived standards are not met, the characteristics of SOP and SPP \\ can make individuals more prone to depression [25,46].} \\ 
 \hline
 \textbf{Dep2} & OOP is tied to lower burnout [7] and positive self-regard [49]. \\ 
 \hline
\end{tabular}

\section{Mathematical Model and Optimization Algorithm}
PSA is initialized with random solution vectors (individuals or persons). For generating new solution, the
algorithm iterates over two main steps called Striving and Depression. Striving step is done in two phases. In phase
one, the PSA selects one of the three types of perfectionism SOP, OOP and SPP according to their associated
probabilities $P_{SOP}, P_{OOP}, P_{SPP}$, respectively. In phase two, a novel solution is produced with regards to the
mathematical model of the selected type of perfectionism. In order to achieve this, SOP, OOP and SPP are
individually mathematically modeled as three operators, based on their characteristics. \\
After generating a new solution, the depression step is executed. If the new solution has a desirable performance
(i.e., it is better than the previously obtained worst solution in the population), algorithm continues to progress by
randomly selecting a type of perfectionism and generating another solution. Otherwise, the probability of the selected
type is updated. Fig. 1 shows the flow chart for a typical iteration performed by PSA. Moreover, in order to further
improve the performance of the algorithm (in terms of complexity), PSA has been made to control the number of
function evaluations by generating as many new solutions as the population size. 

\begin{figure}[ht]
\begin{center}

	\includegraphics[height=5cm]{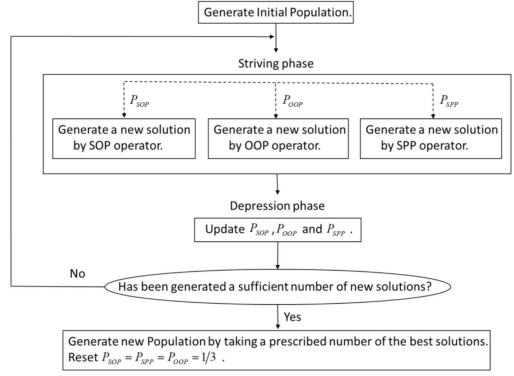}
	\caption{Outline of PSA}
     
\end{center}
\end{figure}

In the following subsections, the details of the different steps, consisting of generating the initial individuals and
performing the striving and depression steps, will be presented. Then, the mathematical models of SOP, OOP and
SPP are provided. The detailed description of how the depression step is created will also be presented later,
followed by the complete description of the PSA algorithm. 

\subsection{Initialization}
As mentioned earlier, the first step is to randomly generate the population over the solution space. The objective function value for each solution is treated as the performance of that solution, and PSA focuses on achieving solutions with high levels of performance. In an n-dimensional optimization problem, a solution (individual or person) is represented as $x = (x_{1},x_{2},...,x_{n})$ where $n$ is the number of the decision variables (such as activities,
traits or behavioral factors of that person). The performance of each solution $x = (x_{1},x_{2},...,x_{n)}$ is computed by evaluating its objective function value, $f(x) = f(x_{1},x_{2},...,x_{n})$. Also, the number of solution vectors in the
population is denoted by $N$ which is kept constant during the optimization process. During the searching process, the population update is accomplished by adding a set of newly generated solutions to the population and then removing the same number of worst solutions from the population; this is done in order to maintain the same number for the total population size ($N$). This process ensures that only the best solutions are kept within the population, therefore adhering to our main inspiration, i.e., perfectionism - as the name suggests. At the end of each iteration, solutions $X_{1}, X_{2} ,..., X_{N}$ are put in order (sorted) from the best to worst according to their objective values (their performance values), i.e., for a minimization problem: $f(X_{1}) \le ... \le f(X_{N})$. So, $X_{1}$ is the best solution (the solution with the smallest objective value), $X_{2}$ is the solution with the second smallest value and so forth, up to X N
that is the solution with the largest objective value (i.e. the worst solution).\\
As mentioned before, for each new solution, the generation process is carried out over two phases. Phase one consists of choosing one of the three types of perfectionism. Phase two consists of using the characteristics from the chosen type to generate a new solution. In the following section, the mathematical formulations for SOP, OOP and SPP as optimization operators are explained in more details.

\subsection{Self-Oriented Perfectionism (SOP)}
SOP operator generates a new solution $X' = \{x'_{1},x'_{2}, ... , x'_{n})$ by modeling the characteristics SOP1 and SOP2, which were defined in section 2. In this case, new solution $X'$ is considered as a self-oriented perfectionist who expects to be perfect ( SOP1) and has exceedingly high personal standards ( SOP2). For this purpose, $10\%$ of the best solutions (i.e., $10\%$ of the first solutions in the sorted population) are considered as a set of patterns for $X'$ The $X'$ perfectionist, in the first instance, randomly selects a solution (person) $X^{*} = (x^{*}_{1}, x^{*}_{2}, ... , x^{*}_{n})$ as a finite sequences $\{(j, x^{*}_{j})\}^{n}_{j=1}$ including $n$ terms in the xy-plane. So, finite sequence $\{(j, x^{*}_{j})\}^{n}_{j=1}$ can be shown as number of points in the plane where the horizontal axis ($x$) is the index number of the term j , and the vertical axis ($y$) is its value, $x^{*}_{j}$ . Since this finite sequence is a function whose domain is set $1,2,...,n$, its graph (denoted by $G^{*}$) consists of some isolated points with coordinates $(1, x^{*}_{1}), (2, x^{*}_{2}), ... , (n, x^{*}_{n})$. In PSA algorithm, Graph $G^{*}$ is interpreted as the behavior of $X^{*}$ . An example of how such a (behavioral) graph may look like is depicted in Fig. 2.

\begin{figure}[ht]
\begin{center}

	\includegraphics[height=5cm]{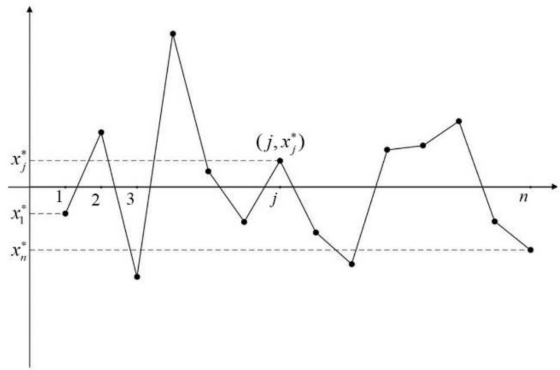}
	\caption{Graph $G^{*}$ : the interpretation of the behavior of point $X^{*}$ .}
     
\end{center}
\end{figure}

A new solution $X'$ is generated by imitating the behavior of $X^{*}$ . In doing so, a random pattern X * is first selected from the set of patterns. Then, SOP operator chooses a random component $j_{0} \in \{1,2,...,n\}$ , and assigns the value of the j‘th component from the selected solution ( $x^{*}$ ), to the corresponding component of the new solution ( $x'_{j_{0}}$ ). Thus, we will have $(j_{0}, x^{*}_{j_{0}}) = (j_{0}, x'_{j_{0}})$. Point $(j_{0}, x^{*}_{j_{0}})$ (or equivalently, $(j_{0}, x'_{j_{0}}$) is labelled by $^{"*"}$ in Fig. 3. Finally, other components $x'_{k} (k \neq j_{0})$ are then generated by the following equation: 

\begin{equation}
    x'_{k} = x'_{j_{0}} + r_{k} (x^{*}_{k} - x^{*}_{j_{0}}), ~~~~ k \in \{1,2,...,n\}-\{j_{0}\}
\end{equation}

Where $r_{k}$ is a uniform random number between $-2$ and $2$ . Algorithm 1, as shown below presents the generation of a new solution $X' = \{x'_{1},x'_{2}, ... , x'_{n})$ using SOP operator. \\

\begin{tabular}{l} \hline \hline
    Algorithm 1 (SOP operator) \\ \hline \hline
    Randomly select solution (pattern) $x^{*}$ from $\{X_{1}, X_{2}, ... , X_{\frac{N}{10}}\}$\\
    Randomly choose a component $x^{*}_{j_{0}} (j_{0}\in \{1,2,...,n\})$ from the selected solution $X^{*} = (x^{*}_{1}, x^{*}_{2}, ... , x^{*}_{n})$. \\
    Set $x'_{j_{0}} = x^{*}_{j_{0}}$\\
    \textbf{for} each $k \in \{1,2,...,n\} - \{j_{0}\}$\\
    ~~~~ Choose a uniform random number; $r_{k} \in U(-2,2)$\\
    ~~~~ Generate $x'_{k}$ byEquation(1).\\
    \textbf{end for}\\
    \hline \hline
\end{tabular}
\\
\\
A clear visualization of how a new position is generated using SOP operator is presented in Fig. 3, where in graph
$G^{*}$ , the behavior of $X^{*}$ is shown by solid line and in graph $G'$ , the behavior of $X'$ is shown by dotted line, which
are all sawtooth-shaped. Graph $G^{*}$ includes points $(k , x^{*}_{k} )$ – shown by black circles – and $G'$ consists of points k
$(k , x'_{k} )$ – shown by white circles.

\begin{figure}[ht]
\begin{center}

	\includegraphics[height=5cm]{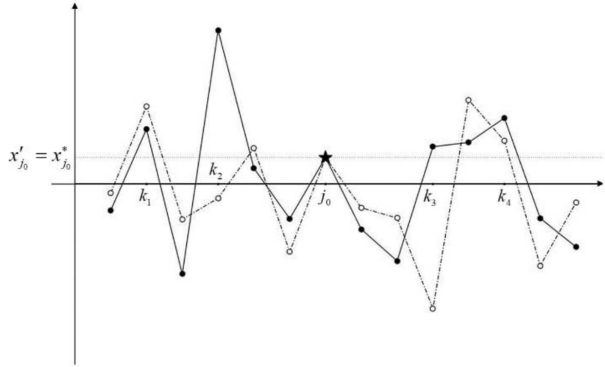}
	\caption{Graph $G'$ : the imitation of the behavior of $X^{*}$ by new solution $X'$ .}
     
\end{center}
\end{figure}

To illustrate the rationale behind equation (1), graphs $G^{*}$ and $G'$ need to be considered, and since $x'_{j_{0}} = x^{*}_{j_{0}}$, two points $(j_{0}, x^{*}_{j_{0}}) \in G^{*}$ and $(j_{0}, x'_{j_{0}}) \in G'$ coincide and are equal to each other. Then, $G'$ is obtained from $G^{*}$ by replacing every point of $(k, x^{*}_{k}) \in G^{*}$ by point $(k, x^{*}_{j_{0}}+r_{k}(x^{*}_{k}-x^{*}_{j_{0}})) $ whose vertical difference from the point $(j_{0}, x^{*}_{j_{0}})$ has been enlarged if $r_{k} > 1$ and reduced if $r_{k} < 1$. This construction can be called $r_{k}$-fold vertical expansion of the graph $G^{*}$, where the term “expansion” allows for enlargement $(r_{k} > 1$, e.g., point $(k_{1}, x^{*}_{k_{1}})$ as shown in Fig.3). no change at all ($r_{k} = 1$, e.g., point $(j_{0}, x^{*}_{j_{0}})$) and reduction ($0<r_{k}<1$, e.g., point $(k_{4}, x^{*}_{k_{4}})$).  Moreover, it also allows for enlargement or reduction coupled with reflection $(r_{k} < 0)$ in the x-axis, that is, $-2 < r_{k} < -1$ for reflected enlargement (e.g., point $(k_{3}, x^{*}_{k_{3}})$) and $-1 < r_{k} < 0$ for reflected reduction (e.g., point $(k_{2}, x^{*}_{k_{2}})$).\\
Consequently, in order to mathematically model the behavior of self-oriented perfectionists, condition SOP1 is achieved by generating a new solution according to one of the best solutions, which is selected from the top 10 percent of solutions. Also, Equation (1) satisfies condition SOP2 in the sense that it can be considered as a high personal standard performed to each component of the new solution.

\subsection{Other-oriented perfectionism (OOP)}
OOP operator generates new solutions based on the best solution obtained so far (the most perfect person) within the population. Similar to previously, the population is sorted at the end of each iteration. So, the existing best solution is
the first solution (i.e., $X_{1} = \{x_{1,1}, x_{1,2}, ..., x_{1,n}\}$) in the sorted population. By this operator, $X_{1}$ is considered as the other-oriented perfectionist. For generating a new solution, OOP operator focuses only on $X_{1}$ (condition OOP4, i.e., narcissism) and does not pay attention to other solutions of the population. Also, it does not allow other solutions $(X_{2}, X_{3}, ... , X_{n})$ to participate directly in the production of new solutions (condition OOP3 , i.e., a lack of trust). The fundamental idea underlying OOP operator is the generation of a new solution around $X_{1} = \{x_{1,1}, x_{1,2}, ..., x_{1,n}\}$ (condition OOP1, where $X_{1}$ expects $X'$ (new solutions) to be perfect). To model the only remaining condition ( OOP2 , i.e., $X_{1}$ has high standards for $X'$ ), some restrictions are created on the position of $X'$ , and as such the
j’th component of a new solution $X' = \{x'_{1}, x'_{2}, ... , x'_{n}\}$ is calculated as follows:
\begin{equation}
    x'_{j} = x_{1.j} + n(0, 2\sigma_{j}), ~~~~ j = 1,2,...,n
\end{equation}
Where $n(0, 2\sigma_{j})$ is a randomly generated number with a normal distribution, mean value of 0 and standard deviation of $2\sigma_{j}$ .  In order to establish the value of standard deviations, the average distance from the j ’th
component of best solution to that of other solutions $X_{1} = \{x_{1,1}, x_{1,2}, ..., x_{1,n}\}$ , $i = 1,2,...,K_{best}$ are calculated, where $K_{best} \in \{2,3,...,N\}$. So, using approach, the distance between $x_{1,j}$ and all $x_{i, j}$ ‘s , $i = 1,2,...,K_{best}$ are initially calculated, and then the value of parameter $\sigma^{2}_{j}$ is created as follows:
\begin{equation}
\frac{\sum_{i=2}^{K_{best}}\vert x_{1,j} - x_{i,j} \vert }{N-1}
\end{equation}
Where $\vert x_{1,j} - x_{i,j} \vert$ indicates the distance from j’th component of $X_{1}$ (the best solution obtained so far) to that of $X_{i}$ . The whole process is repeated for each dimension $j = 1,2,...,n$ and each time the average distance is
calculated using only the single dimension j . An example of how this approach can be used to generate the j‘th
component $(x'_{j})$ from the new solution $X' = \{x'_{1}, x'_{2}, ... , x'_{n}\}$, is presented in Fig. 4. This figure illustrates equation (2), where $P(x'_{j})$ is the probability that value $x'_{j}$ will be selected as the j’th component of $X'$.

\begin{figure}[ht]
\begin{center}

	\includegraphics[height=5cm]{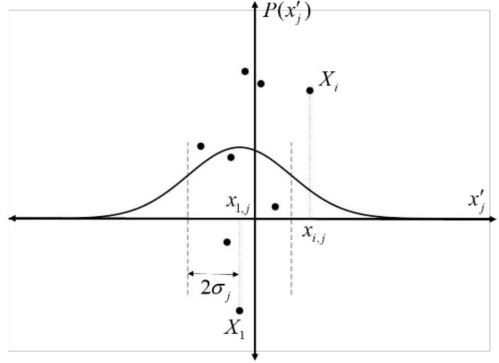}
	\caption{Probability density function $P(x'_{j})$ for generating j’th component of $X'$.}
     
\end{center}
\end{figure}

The position of $X_{1}$ and n–dimensional cube generated by $\sigma_{j}$‘s $(j = 1,2,...,n)$, for $n=2$ has been shown in Fig. 5. According to equations (2) and (3), the probability that a certain point $X' = \{x'_{1}, x'_{2}, ... , x'_{n}\}$ is generated as a new solution is equal to:
\begin{equation} 
\Pi_{j=1}^{n}\{ \frac{1}{2\sigma_{j}\sqrt{2\pi}} e^{-\frac{1}{2}(\frac{x'_{j}-x_{1,j}}{2\sigma_{j}})^{2}} \}
\end{equation}
Fig. 5 shows these probabilities for all the points in the plane. The points located in the darker regions are more likely to be selected as the new solutions.
\begin{figure}[ht]
\begin{center}

	\includegraphics[height=5cm]{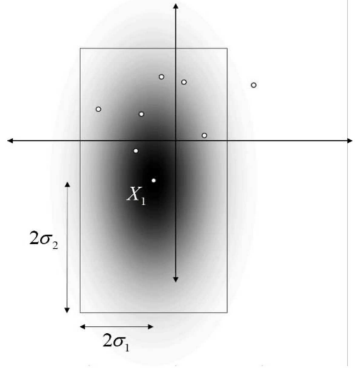}
	\caption{Cube generated by $\sigma_{1}$ and $\sigma_{2}$ ; the probabilities of points which may be generated as a new solution.
}
     
\end{center}
\end{figure}

During the entire process of iterations, new solutions which are generated by the OOP operator, explore a good area surrounding the best solution in a random manner. However, using the standard normal distribution presented in equation (2), the search is more biased toward the best solution and less biased towards points which their distances from the best solution are greater than the average distance calculated by equation (3). This guarantees the exploitation step of the algorithm. Summary of how Algorithm 2 (OOP operator) can be expressed is as follows:

\begin{tabular}{l} \hline \hline
    Algorithm 2 (OOP operator) \\ \hline \hline
    Choose a uniform integer random number $K_{best} \in \{1,2,...,N\}$\\ 
    \textbf{for} $j = 1,2,...N$\\
    ~~~~ Calculate $\sigma_{j}$ by Equation (3).\\
    ~~~~ Generate a random number $x'_{j}$ by Equation (2)..\\
    \textbf{end for}\\
    \hline \hline
\end{tabular}
\\
\\

The stochastic character of variable $K_{best}$ can be used as the way to achieve a reasonable compromise between exploration and exploitation steps. In such a way, only the best solutions ( $K_{best}$ ) are used in equation (3). The lower the value of $K_{best}$, the higher the convergence speed of the algorithm. In this case, the worse solutions in the population are not considered.

\subsection{Socially Prescribed Perfectionism (SPP)}
On the contrary to the previous approach, the SPP operator generates a new solution corresponding to the solutions from the population. In this case, a new solution $X'$ is interpreted as a socially prescribed perfectionist who is influenced by all the solutions of the population. The main idea behind SPP consists of two steps: (1) all the solutions have the same chance to generate $X'$ (condition SPP1, i.e., $X'$ believes others expect him/her to be perfect); (2) the components of $X'$ are the copies of the corresponding components from the solutions $X_{1}, X_{2}, ... , X_{N}$ (condition SPP2 , i.e., $X'$ desires approval from those around him/her). Therefore, SPP operator is simply stated as follows:\\

\begin{tabular}{l} \hline \hline

    Algorithm 3 (SPP operator) \\ \hline \hline
    \textbf{for} $j = 1,2,...N$\\
    ~~~~ Randomly select a solution $\bar{X} = \{ \bar{x_{1}}, \bar{x_{2}}, ... , \bar{x_{N}}  $.  \\
    ~~~~ Set $x'_{j} = \bar{x_{j}}$.\\
    \textbf{end for}\\
    \hline \hline
\end{tabular}
\\
\\
This operator generates the new solutions within an area of minimal cell ( n-dimensional cube) including all the previously obtained solutions $(X_{1}, X_{2}, ... , X_{N})$. So, this operator enhances the diversity of the population and shares information among individuals. As a consequence, it assists the proposed algorithm to search solution space randomly and avoid any trap within the local optima.

\subsection{Depression}

The fundamental idea underlying the depression step is the formulation of depression and its effects on the self- oriented and socially prescribed perfectionists (see condition Dep1). Depression can be mathematically modelled based on the biased probabilistic choice of SOP and SPP operators. Furthermore, due to the condition Dep2 , no negative effects have been considered for OOP operator in this step. These probabilities are initialized with the initial values of $P_{SOP} = P_{OOP} = P_{SPP} = \frac{1}{3}$. During an arbitrary iteration, if a new solution ($X'$) which is generated by SOP operator fails to meet the set expectation (i.e., $X'$ is worse than the worst solution within the current population), then probability ($P_{SOP}$) is decreased by a certain value of $\Delta = \frac{1}{3N}$. To make sure that the equation
$P_{SOP} + P_{OOP} + P_{SPP} = 1$ remains true for the next generation, the probabilities of $P_{SPP}$ and $P_{OOP}$ are increased to $P_{SPP} + \frac{\Delta}{2}$ and $P_{OOP} + \frac{\Delta}{2}$, respectively. So, for generating new solution, SOP , OOP and SPP may be selected with probabilities of $P_{SOP} - \Delta$, $P_{OOP} + \frac{\Delta}{2}$, and $P_{SPP} + \frac{\Delta}{2}$, respectively. Similarly, for each of the failed new solution ($X'$) generated by SPP operator, $P_{SPP}$ is decreased to $P_{SPP} - \Delta$; $P_{SOP}$ and $P_{OOP}$ are increased by $\frac{\Delta}{2}$. This process will continue until the end of the iteration.

As mentioned earlier, N number of new solutions are produced during each iteration and, therefore, each of the SOP, SPP and OOP operators can be selected at most N times during each iteration. Hence, in an arbitrary iteration, if SOP (or SPP) operator is selected N times with N number of failed new solutions, then the associated probability will be reduced from the initial value $\frac{1}{3}$ to $\frac{1}{3} - N\Delta = 0$ at the end of the iteration. However, at the beginning of the next iteration, we set $P_{SOP} + P_{OOP} + P_{SPP} = \frac{1}{3}$. Using this strategy, the algorithm will have the same opportunity (in each iteration) to search the space using the three operators. Hence, Depression step makes the algorithm more likely to choose the most successful types of operators.\\
The probability update, however, is optional in the sense that $P_{SOP}$,  $P_{OOP}$, and $P_{SPP}$ may be considered as the same constant values over the course of iterations. In this case, any of the perfectionism type has the same chance of being selected, i.e., each of the types (SOP, OOP, or SPP) has a probability of $33.3 \%$ being selected. On the one hand, this strategy can enable the algorithm to quickly compute approximate solutions, but on the other hand, it causes the algorithm to achieve less accurate solutions. In the second way, the probabilities can be updated during an iteration. Although this method allows carefully considered selections of SOP, OOP and SPP operators (and hence increases the), but it also decreases the speed of the algorithm. In this paper, we have employed a method to compromise between the accuracy of the algorithm and the speed at which the best solutions can be achieved. To this end, we set $P_{OOP} = \frac{1}{2}$ and $P_{SOP} = P_{SPP} = \frac{1}{4}$ as the fixed probabilities for all iterations.\\
Finally, Algorithm 4 shows the PSA algorithm for a minimization problem, where $OP \in \{SOP,OOP,SPP\}$ is a variable indicating which operator has been selected for generating a new solution.\\

\begin{tabular}{l} \hline \hline
    Algorithm 4 (Perfectionism Search Algorithm) \\ \hline \hline
    Generate the initial population.\\
    Evaluate each solution and sort the solutions in an ascending order of their objective values.\\
    Calculate the best solution $X^{*}$ by setting $X^{*} = X_{1}$.\\
    \textbf{While} ( $t \le$ maximum number of iterations) \\
    ~~~~\textbf{for} $i = 1,2,...N$\\
    ~~~~~~~~ Select operator $OP \in \{SOP,OOP,SPP\}$ according to probabilities \\
    ~~~~~~~~ $P_{OOP} = \frac{1}{2}$ and $P_{SOP} = P_{SPP} = \frac{1}{4}$.\\
    ~~~~~~~~~~~~ \textbf{if} OP = SOP, \\
    ~~~~~~~~~~~~~~~~ Generate new solution $X'$ by Algorithms 1 and evaluate $X'$.\\
    ~~~~~~~~~~~~ \textbf{elseif} OP = OOP, \\
    ~~~~~~~~~~~~~~~~ Generate new solution $X'$ by Algorithms 2 and evaluate $X'$.\\
    ~~~~~~~~~~~~ \textbf{else}\\
    ~~~~~~~~~~~~~~~~ Generate new solution $X'$ by Algorithms 3 and evaluate $X'$.\\
    ~~~~~~~~~~~~ \textbf{end if}\\
    
    ~~~~ \textbf{end for}\\
    ~~~~Generate the new population by keeping \\ ~~~~the N best solutions and removing the N worst ones.\\
    ~~~~Sort the solutions of the resulted population.\\
    ~~~~Update the best solution $X^{*}$ by setting $X^{*} = X_{1}$ , if $f(x_{1}) \le f(x^{*})$. \\
    ~~~~$t = t + 1$. \\
    \textbf{end while} \\
\hline \hline    
\end{tabular}
\\
\\
\section{Theoretical Aspects, Applied Suggestions and Brief Comparative Study}

In order to see how PSA solves optimization problems, it’s important to note some points:
\begin{enumerate}[label=(\Alph*)]
    \item In PSA, new solutions are generated using the three types of optimization operators. It assists the algorithm
    to look for an optimal point from different strategies.
    \item In PSA, each search strategy focuses on a specific region; SOP considers neighborhoods around $10 \%$ of the
    best solutions, OOP pays attention to the neighborhood of the best-so-far solution and SPP includes the region between the entire solutions. As such, it creates a balance between exploration and exploitation, and increases the capability of the algorithm to be applied to the optimization of multi-modal problems.
    \item By considering $P_{OOP} = \frac{1}{2}$ and $P_{SOP} = P_{SPP} = \frac{1}{4}$, an appropriate balance is made between the accuracy and speed of the algorithm.
    \item In SOP operator, if $\vert r \vert$ is close to 2 ($\vert r \vert  \approx  2$), the exploration power is increased around the current best so-far solution, which allows for more diversification of the search (higher robustness). Also, if $\vert r \vert$ is
    close to 0 ($\vert r \vert  \approx  0$), the intensification capability is increased, that leads to a faster convergence (higher
    convergence speed and higher efficiency)
    \item In OOP, since new solutions are generated from the neighborhoods of $10 \%$ of the best solutions, more
    promising regions are investigated.
    \item SOP and OOP assist PSA to exploit information from the neighborhood of the best solutions that hold valuable knowledge.
    \item In SOP (when $\vert r \vert < 1$) and OOP (by lapse of iterations), the space is searched more locally. So, it can be
    considered as an adaptive learning rate.
    \item The minimal covering cell used by SPP operator assists the algorithm to increase the exploration power. So, the search mechanism in SPP has two advantageous; escaping from local optima and avoiding rapid convergence of the algorithm.
    \item Since in SPP, all the previously obtained solutions contribute towards generation of new solutions, this operation can be considered as an information – transferring tool.
\end{enumerate}
The following points highlight some of the main characteristics of PSA that differentiates it from the other algorithms.
\begin{enumerate}[label=(\roman*)]
    \item Unlike many of the previously introduced heuristic algorithms (such as ACO, PSO, CA, HS, etc), PSA is a memory-less algorithm that works efficiently and has a good convergence rate.
    \item One major difference between PSA and most of the other heuristic algorithms is that many heuristic algorithms use some improving directions (e.g., PSO, BBO, GSA, FA, BA, TLBO, BH, ...), while PSA does not apply any particular direction to improve the solutions.
    \item SOP operator, unlike many optimization operators, does not consider the distance between solutions; however, it works only by changing the differences between the components of the best solution.
\end{enumerate}
\section{Experimental Results}
In order to evaluate the performance of our algorithm in comparison to other known methods, the typical benchmark test functions have been employed by selecting a comprehensive set of 35 benchmark functions that are typically used in the literature for presenting the performance of different methods and algorithms. List of the test functions, number of dimensions used, initialization interval, optimal solutions and the optimal values are presented in Appendix A. A more detailed description of the test functions may be found in [28,63]. These functions listed in Appendix, Tables A.1-A.3 range from very simple to quite complex, and were selected for their particular characteristics. Functions $f_{i} (i = 1,2,...,11)$ are unimodal (Table A.1) and $g_{i} (1 = 1,2,...,13)$ are multimodal (Table A.2). Also, $\widetilde{h_{i}} (i = 1,2,...,11)$ are the shifted functions which have resulted from some unimodal and multimodal functions $h_{i}$ so that $\widetilde{h_{i}}(x) = h_{i}(x - \textbf{1})$ where $\textbf{1} = [1,1,1,...,1]_{1 \times n}$.

The proposed PSA is compared with 11 well-known metaheuristic methods including the most famous algorithms and also some of the recent popular ones, namely, Ant Colony Optimization (ACO) [6,48], Particle Swarm Optimization (PSO) [29], Teaching Learning Based Optimization (TLBO) [40], Biogeography Based Optimization (BBO) [44], Differential Evolution (DE) [53], Harmony Search (HS) [16], Firefly Algorithm (FA) [61], Whale Optimization Algorithm (WOA) [35], Bees Algorithm (BA) [38], Gravitational Search Algorithm (GSA) [41] and Black Hole (BH) algorithm [22]. In order to ensure a fair comparison, we carefully used the parameters suggested by the authors for each algorithm.

In all cases, population size (N) is set to $N = 50$ . 30-dimensional versions of the functions ( $n = 30$ ) were used and 30 independent runs were performed in all the experiments. In addition, two types of comparison were made. For the first type, maximum iteration is considered as a stopping condition and is set to 1000 . For each test problem, the following performance indices (averaged over 30 runs) are reported: the best optimal value (Table 1), the average best-so-far solutions (Avg), median of the best solutions (Mdn) and standard deviation (SD) for the last iterations (Tables 2 – 4). For the second type, a comparison was made based on the number of function evaluations needed to reach the following stopping condition:
\begin{equation}
    \vert f_{best} - f^{*} \vert < \epsilon f_{best} + \epsilon_{2}
\end{equation}
Where $f_{best}$ is the value of the best solution found by an algorithm, $f^{*}$ is the known optimal value for a given test problem,
and $\epsilon_{1}$ and $\epsilon_{2}$ are respectively the relative and absolute errors [48]. For all the test runs, we used $\epsilon_{1} = \epsilon_{2} = 10^{-4}$. For each test problem, the mean number of function evaluations (NFE) is supplied for all the algorithms. Moreover, numbers in square brackets indicate the percentage of successful runs (PSR), i.e., average number of cases in which the required accuracy was reached (Tables 2 – 4). In all the tables, shaded boxes have been used to show the best results for a given problem.

\hspace{-2cm}
\begin{table}[H]
    \centering
    \begin{tabular}{llllllll}
\hline \hline
 &         WOA &        TLBO &        BBO &         DE &         ACO &          PSA \\
\hline \hline
$f_{1}$ & 4.7485e-276 & 3.7512e-246 & 1.0295e-10 & 7.9684e-19 &  6.6077e-72 &            0 \\ \hline
$f_{2}$ &      0.6667 &      0.9041 &     0.3490 &     0.6667 &      0.6667 &   2.2288e-19 \\\hline
$f_{3}$ & 3.6732e-266 & 9.8671e-245 & 6.1491e-08 & 6.9791e-10 &  1.3055e-09 &            0 \\\hline
{$f_{4}$} &         0.5 &         0.5 &        0.5 &        0.5 &         0.5 &          0.5 \\\hline
{$f_{5}$} &     26.2707 &     28.7923 &     5.2982 &    18.3086 &     0.71274 &            0 \\\hline
{$f_{6}$} &     34.2804 & 1.5839e-243 &  1.5839e02 &    0.98238 &   0.0011337 &            0 \\\hline
{$f_{7}$} & 5.7295e-149 & 6.2017e-121 & 2.1409e-05 & 3.6638e-08 &  4.8992e-41 &  3.2188e-245 \\\hline
{$f_{8}$} &  3.402e-152 & 1.3215e-121 & 2.1443e-06 &  7.268e-09 &  5.6872e-42 &  5.8602e-247 \\\hline
{$f_{9}$} & 1.6524e-273 & 7.6983e-245 & 4.3522e-09 &   8706e-15 &   6.416e-69 &            0 \\\hline
{$f_{10}$} &           0 &           1 &          0 &          0 &           0 &            0 \\\hline
{$f_{11}$} & 2.0678e-275 & 3.0112e-244 & 3.8436e-10 & 1.1107e-16 &  4.3393e-71 &            0 \\\hline
{$g_{1}$} &  8.8818e-16 &  4.4409e-15 & 1.1597e-05 & 1.2372e-08 &  4.4409e-15 &   8.8818e-16 \\\hline
{$g_{2}$} &    -0.86175 &     -0.7625 &         -1 &   -0.76236 &    -0.66986 &           -1 \\\hline
{$g_{3}$} &          -1 &    -0.93625 &   -0.47778 &   -0.81731 &    -0.78575 &           -1 \\\hline
{$g_{4}$} &           0 &           0 & 3.6267e-08 & 2.6645e-15 &           0 &            0 \\\hline
{$g_{5}$} &      0.2571 &      1.3853 & 5.2297e-11 & 2.4151e-16 &      0.0895 &   1.4998e-32 \\\hline
{$g_{6}$} &           0 &  2.9848e-07 &     4.0758 &          0 &      11.707 &            0 \\\hline
{$g_{7}$} &      0.2409 &      1.7004 & 4.3376e-10 & 1.4941e-15 &  2.3359e-32 &   1.3498e-32 \\\hline
{$g_{8}$} &  9.146e-275 & 2.3084e-243 &     0.1475 & 7.3163e-09 &   2.795e-51 &            0 \\\hline
{$g_{9}$} &      12.533 &     47.1373 & 3.7756e-05 & 4.8925e-06 &  4.4409e-16 &            0 \\\hline
{$g_{10}$} &  3.618e-134 &      0.0999 &     0.5999 &     0.1999 &      0.1999 &            0 \\\hline
{$g_{11}$} & 2.0415e-146 & 8.7349e-243 &     0.0001 & 1.5972e-11 &  4.7422e-50 & -1.2302e-320 \\\hline
{$g_{12}$} &     -3.4551 &      -3.165 &       -3.5 &    -1.4239 &     -0.0657 &         -3.5 \\\hline
{$g_{13}$} &    -38.8512 &    -27.1919 &   -37.7525 &   -39.1662 &    -37.2813 &     -39.1662 \\\hline
{$\widetilde{h_{1}}$} &      0.0942 &      5.9096 & 4.5125e-07 &  1.258e-07 &  3.0705e-16 &            0 \\\hline
{$\widetilde{h_{2}}$} &   1.2157e05 &   1.1110e07 &     0.0028 & 1.3957e-09 &  3.6978e-26 &            0 \\\hline
{$\widetilde{h_{3}}$} &     22.2976 &     50.3096 & 5.0246e-08 & 6.8945e-14 &           0 &            0 \\\hline
{$\widetilde{h_{4}}$} &      0.1101 &   1.7681e02 & 4.6465e-17 & 6.5906e-30 &  5.4694e-63 &            0 \\\hline
{$\widetilde{h_{5}}$} &  1.4442e-05 &      6.9221 &  2.889e-20 & 4.3705e-14 &  2.2613e-88 &   1.6005e-94 \\\hline
{$\widetilde{h_{6}}$} &   1.4569e02 &   3.6972e03 &     0.0004 & 1.6213e-13 &  5.1465e-30 &            0 \\\hline
{$\widetilde{h_{7}}$} &     26.0988 &     29.9126 &    13.9294 & 10.6000e01 &     19.8992 &            0 \\\hline
{$\widetilde{h_{8}}$} &  6.9944e-05 &      5.4866 & 4.3223e-14 & 8.8788e-32 & 2.5627e-144 &            0 \\\hline
{$\widetilde{h_{9}}$} &      5.6469 &     33.3308 &     2.0375 &      0.128 &      0.0023 &            0 \\\hline
{$\widetilde{h_{10}}$} &      0.3005 &      0.5197 &     0.2422 &     0.4975 &     0.06738 &            0 \\\hline
{$\widetilde{h_{11}}$} &   3.7697e02 &     19.4955 &    10.8001 &     0.1280 &  1.0532e-06 &            0 \\\hline \hline

\end{tabular}
    \caption{THE BEST OPTIMAL VALUES OF THE TEST PROBLEMS FOUND BY WOA, TLBO, BBO, DE, ACO AND PSA}
\end{table}

\begin{table}[H]
    \centering
\begin{tabular}{llllllll}
\hline \hline
{} &         FA &         BA &        PSO &         HS &         BH &        GSA \\ \hline \hline

{$f_{1}$} &     2.3386 & 2.2448e-23 & 2.6202e-16 &     0.0112 &     1.6336 &  4.4408e02 \\\hline
{$f_{2}$} &    26.7541 &     0.6667 &     0.6668 &     6.4017 &     0.6667 &     0.6667 \\\hline
{$f_{3}$} &  1.4256e02 &     0.0012 & 8.8712e-07 &     2.1427 & 6.5624e-05 & 2.6722e-06 \\\hline
{$f_{4}$} &     3.1697 &      2.353 &        0.5 &     0.5000 &        0.5 &        0.5 \\\hline
{$f_{5}$} & 18.1205e01 &    19.6798 &    25.8392 & 24.7479e01 &     26.766 &    25.1698 \\\hline
{$f_{6}$} &     5.4192 &    17.4221 & 91.5358e02 & 22.0260e02 &    86.5808 &    28.7561 \\\hline
{$f_{7}$} &     5.8056 &    16.0956 & 1.9953e-07 &    10.3322 & 4.8536e-06 & 1.0754e-08 \\\hline
{$f_{8}$} &     5.2023 &     2.9572 & 1.0309e-08 &    0.79269 & 3.3232e-08 & 1.0185e-08 \\\hline
{$f_{9}$} &     1.5422 & 1.9182e-20 & 2.2316e-15 &     9.0222 & 5.7282e-09 & 5.1986e-18 \\\hline
{$f_{10}$} &          0 &         25 &          0 &         10 &          0 &          0 \\\hline
{$f_{11}$} &    19.8964 & 1.0639e-08 & 2.4254e-16 &     1.1945 & 1.0761e-09 & 2.9106e-17 \\\hline
{$g_{1}$} &    18.9259 &    10.5707 &     0.0045 &     1.1693 & 1.5185e-05 & 1.8428e-09 \\\hline
{$g_{2}$} &    -0.7042 &         -1 &   -0.97548 &   -0.98858 &   -0.99636 &         -1 \\\hline
{$g_{3}$} &    -0.6195 &   -0.15277 &   -0.17987 &   -0.78568 &    -0.6195 &   -0.47773 \\\hline
{$g_{4}$} &  1.4024e02 & 2.3315e-15 & 2.2502e-06 &     1.0714 & 1.0145e-06 &     1.1449 \\\hline
{$g_{5}$} &    18.1618 &    22.9755 & 6.9424e-12 &     0.0437 & 6.1388e-10 & 2.0405e-18 \\\hline
{$g_{6}$} &    10.6185 &    10.0504 &    10.3806 &     5.0166 &     9.0053 &     7.7226 \\\hline
{$g_{7}$} &     0.2776 &     0.0141 &      4.594 &     0.6096 & 9.0718e-08 & 5.8712e-19 \\\hline
{$g_{8}$} & 13.6073e02 & 44.2623e02 & 51.2685e01 & 10.2741e01 & 43.7389e01 &    20.9931 \\\hline
{$g_{9}$} &    36.5972 &     3.5685 &    10.0071 &     10.821 & 1.1311e-05 &   1.18e-07 \\\hline
{$g_{10}$} &    0.20869 &    0.59987 &     2.4009 &    0.96516 &    0.49987 &    0.90598 \\\hline
{$g_{11}$} &    32.4466 & 24.9864e04 &     0.0002 & 36.1400e01 & 5.4209e-05 & 9.3929e-17 \\\hline
{$g_{12}$} &    -1.1778 &    -3.4994 &  -0.051724 &   -0.64447 &    -3.4765 &       -3.5 \\\hline
{$g_{13}$} &   -34.5706 &   -37.2811 &   -35.3964 &    -39.155 &   -36.7842 &   -36.3388 \\\hline
{$\widetilde{h_{1}}$} &    11.8138 &     2.3454 &     0.0008 &     0.1004 &     0.0003 & 1.1154e-09 \\\hline
{$\widetilde{h_{2}}$} &  1.5635e06 &     3.6084 & 7.7902e-08 &  6.9181e06 &     0.0618 &     0.3074 \\\hline
{$\widetilde{h_{3}}$} &    20.4495 &    11.5618 & 2.8914e-07 &    19.2699 &  5.544e-07 &          0 \\\hline
{$\widetilde{h_{4}}$} &      2.738 & 2.8102e-31 &   1.11e-34 &  1.0473e02 &  8.234e-15 & 2.4956e-35 \\\hline
{$\widetilde{h_{5}}$} &     0.0555 & 4.3294e-42 & 1.4731e-08 & 4.2817e-09 &     0.3111 & 5.5874e-54 \\\hline
{$\widetilde{h_{6}}$} &  1.7006e04 &  2.5809e04 & 4.6137e-14 &  3.5254e03 &     24.696 &    73.5399 \\\hline
{$\widetilde{h_{7}}$} & 14.8043e01 & 10.3606e01 &    31.8505 &    10.0787 &    24.2123 &    11.9395 \\\hline
{$\widetilde{h_{8}}$} &     0.0008 & 7.3501e-08 & 1.3389e-16 & 8.7973e-11 &     0.0490 & 1.0409e-15 \\\hline
{$\widetilde{h_{9}}$} &  1.4600e02 &     0.5661 &     4.5534 &     3.1674 &     3.2997 &     0.4865 \\\hline
{$\widetilde{h_{10}}$} &     0.5458 &     0.4065 &     0.1382 &     0.0174 &     0.2365 &     0.6184 \\\hline
{$\widetilde{h_{11}}$} &     1.7421 & 2.4295e-09 &  2.2660e02 &    37.1223 &     0.0412 &    52.5817 \\\hline \hline 

\end{tabular}

    \caption{The best optimal values of the test problems found by FA, BA, PSO, HS, BH and GSA}
\end{table}

\textbf{Unimodal functions.} Table 1 include the best values,$f^{*}(x)$, found by the algorithms over 30 runs. These tables demonstrate the impressive ability of the proposed algorithm to detect the optimal solutions from the test problems. PSA is a clear winner for
$f_{1} - f_{3}$ , $f_{5} - f_{9}$ and $f_{11}$ . For $f_{4}$ and $f_{10}$ almost all the algorithms found the global optimum.\\
The good convergence rate of PSA can be observed from Table 2. According to these Tables, PSA produced the best results in
terms of Avg, Mdn and SD for all the test problems, with the only exception of $f_{2}$ where ACO produced a better standard deviation. For $f_{1} - f_{3}$ , $f_{5} - f_{9}$ and $f_{11}$ , PSA is the only algorithm that obtained the best Avg and Mdn. For $f_{4}$ , almost all the algorithms performed similarly well and for $f_{10}$ , some algorithms attained the best results. Based on these results, PSA tends to find the global optimum much faster in comparison to other algorithms and hence has a higher convergence rate. Furthermore, it can produce high quality solutions with low average error (good accuracy) and low standard deviation (good stability) for all the unimodal functions.\\
Moreover, using PSA the required accuracy is reached in every run (i.e., $PSR =100\%$) for each unimodal function; and in all cases it has the minimum number of function evaluations. Also, it is to be noted that for  $f_{2}$ and  $f_{5}$ , PSA is the only algorithm
that reached this accuracy.

\begin{table}[H]
\tiny
    \centering
\begin{tabular}{lllllllll}
\hline \hline
 &  &   WOA & TLBO   & BBO & DE & ACO & PSA \\ \hline \hline
  $f_{1}$ & \makecell{Avg \\ Mdn \\ SD \\ NFE [PSR]} & \makecell{1.4377e-257 \\ 1.1616e-266 \\  0 \\ 12,312.3 \textbf{[100\%]}} & \makecell{3.7929e-242 \\ 6.7492e-244 \\  0 \\ 2,601.3 \textbf{[100\%]}} & \makecell{8.0804e-09 \\ 1.4835e-09 \\ 1.6572e-08 \\ 22,622.6 \textbf{[100\%]}} & \makecell{5.4333e-17 \\ 2.9971e-17 \\ 5.5677e-17 \\ 15,740.7 \textbf{[100\%]}} & \makecell{7.7876e-69 \\ 1.7912e-69 \\ 1.584e-68 \\ 4,954.9 \textbf{[100\%]}|} & \makecell{0 \\ 0 \\  0 \\ 600.6   \textbf{[100\%]}} \\ \hline
  $f_{2}$ & \makecell{Avg \\ Mdn \\ SD \\ NFE [PSR]} & \makecell{0.6667 \\ 0.6667 \\ 2.7959e-05 \\ 50,050 [0\%]} & \makecell{0.9550 \\ 0.9605 \\ 0.0210 \\ 100,050 [0\%]} & \makecell{1.209 \\ 0.6903 \\ 0.9522 \\ 50,050 [0\%]} & \makecell{0.7605 \\ 0.6667 \\ 0.4956 \\ 50,050 [0\%]} & \makecell{0.6667 \\ 0.6667 \\ 1.4189e-16 \\ 50,050 [0\%]} & \makecell{3.6769e-13 \\ 9.9902e-17 \\ 1.799e-12 \\ 11,186.2 \textbf{[100\%]}} \\ \hline
  $f_{3}$ & \makecell{Avg \\ Mdn \\ SD \\ NFE [PSR]} & \makecell{6.3149e-27 \\ 6.6672e-188 \\ 3.3934e-26 \\ 13,963.9 \textbf{[100\%]}} & \makecell{4.5583e-240 \\ 6.4553e-242 \\ 0 3,551.8 \textbf{[100\%]}} & \makecell{0.0001 \\ 8.74e-05 \\ 9.1247e-05 \\ 49,949.9 [53.3\%]} & \makecell{9.2893e-07 \\ 9.9539e-08 \\ 1.8022e-06 \\ 22,022 \textbf{[100\%]}} & \makecell{5.3943e-09 \\ 4.8139e-09 \\ 3.0165e-09 \\ 6,631.6 \textbf{[100\%]}} & \makecell{0 \\ 0 \\ 0 \\ 1,051 \textbf{[100\%]}} \\ \hline
  $f_{4}$ & \makecell{Avg \\ Mdn \\ SD \\ NFE [PSR]} & \makecell{0.5 \\ 0.5 \\ 3.3307e-17 \\ 11,186.2 \textbf{[100\%]}} & \makecell{0.5 \\ 0.5 \\ 0 \\ 1,400.7 \textbf{[100\%]}} & \makecell{0.5 \\ 0.5 \\ 8.2834e-17 \\ 5,280.3 \textbf{[100\%]}} & \makecell{0.5 \\ 0.5 \\ 2.7481e-10 \\ 11,436.4 \textbf{[100\%]}} & \makecell{0.5 \\ 0.5 \\ 5.2336e-17 \\ 2,952.9 \textbf{[100\%]}} & \makecell{0.5 \\ 0.5 \\ 0 \\ 325.3 \textbf{[100\%]}} \\ \hline
  $f_{5}$ & \makecell{Avg \\ Mdn \\ SD \\ NFE [PSR]} & \makecell{27.3042 \\ 27.073 \\ 0.7647 \\ 50,050 [0\%]} & \makecell{28.8927 \\ 28.8999 \\ 0.0433 \\ 100,050 [0\%]} & \makecell{64.541 \\ 77.3741 \\ 33.8184 \\ 50,050 [0\%]} & \makecell{37.4842 \\ 25.8895 \\ 23.7393 \\ 50,050 [0\%]} & \makecell{48.8164 \\ 22.0362 \\ 52.6671 \\ 50,050 [0\%]} & \makecell{\textbf{6.6429e-28} \\ 0 \\ \textbf{3.1271e-27} \\ \textbf{4,254.2} \textbf{[100\%]}} \\ \hline
  $f_{6}$ & \makecell{Avg \\ Mdn \\ SD \\ NFE [PSR]} & \makecell{5.5902e03 \\ 4.6003e03 \\ 4.7965e03 \\ 50,050 [0\%]} & \makecell{2.7584e-239 \\ 1.0467e-240 \\ 0 \\ 3,801.9 \textbf{[100\%]}} & \makecell{7.4737e02 \\ 6.4937e02 \\ 3.2609e02 \\ 50,050 [0\%]} & \makecell{9.8156 \\ 6.1309 \\ 9.9271 \\ 50,050 [0\%]} & \makecell{0.0186 \\ 0.0132 \\ 0.0143 \\ 50,050 [0\%]} & \makecell{3.3907e-294 \\ 1.6304e-321 \\ 0 \\ 1,876.9 \textbf{[100\%]}} \\ \hline
  $f_{7}$ & \makecell{Avg \\ Mdn \\ SD \\ NFE [PSR]} & \makecell{3.1658e-142 \\ 2.0906e-144 \\ 1.2868e-141 \\ 14,639.6 \textbf{[100\%]}} & \makecell{6.6069e-119 \\ 1.2272e-119 \\ 1.1014e-118 \\ 6,003 \textbf{[100\%]}} & \makecell{0.0002 \\ 0.0001 \\ 0.0002 \\ 47,072 [36.7\%]} & \makecell{1.768e-07 \\ 1.4899e-07 \\ 1.0872e-07 \\ 35,435.4 \textbf{[100\%]}} & \makecell{1.4602e-39 \\ 6.7007e-40 \\ 1.6051e-39 \\ 8,358.3 \textbf{[100\%]}} & \makecell{8.52e-225 \\ 1.5417e-233 \\ 0 \\ 1,351.3 \textbf{[100\%]}} \\ \hline
    $f_{8}$ & \makecell{Avg \\ Mdn \\ SD \\ NFE [PSR]} & \makecell{2.1075e-141 \\ 7.5248e-147 \\ 1.1259e-140 \\ 14,664.6 \textbf{[100\%]}} & \makecell{1.1094e-119 \\ 1.7977e-120 \\ 2.34e-119 \\ 5,252.6 \textbf{[100\%]}} & \makecell{1.6921e-05 \\ 1.2936e-05 \\ 1.4481e-05 \\ 35,160.1 \textbf{[100\%]}} & \makecell{4.5861e-08 \\ 4.2129e-08 \\ 2.4938e-08 \\ 31,406.4 \textbf{[100\%]}} & \makecell{2.4516e-40 \\ 1.4192e-40 \\ 2.5119e-40 \\ 7,757.7 \textbf{[100\%]}} & \makecell{3.9638e-225 \\ 2.4596e-238 \\ 0 \\ 1,226.2 \textbf{[100\%]}} \\ \hline
  $f_{9}$ & \makecell{Avg \\ Mdn \\ SD \\ NFE [PSR]} & \makecell{1.012e-255 \\ 2.4036e-262 \\ 0 \\ 14,589.6 \textbf{[100\%]}} & \makecell{7.596e-240 \\ 1.039e-241 \\ 0 3,751.9 \textbf{[100\%]}} & \makecell{2.4358e-07 \\ 8.8815e-08 \\ 3.5563e-07 \\ 24,699.7 \textbf{[100\%]}} & \makecell{2.2608e-14 \\ 1.2951e-14 \\ 2.1513e-14 \\ 22,722.7 \textbf{[100\%]}} & \makecell{1.398e-66 \\ 3.0118e-67 \\ 3.0587e-66 \\ 6,281.3 \textbf{[100\%]}} & \makecell{0 \\ 0 \\ 0 \\ 825.8 \textbf{[100\%]}} \\ \hline
  $f_{10}$ & \makecell{Avg \\ Mdn \\ SD \\ NFE [PSR]} & \makecell{0 \\ 0 \\ 0 \\ 10,310.3 \textbf{[100\%]}} & \makecell{4.7667 \\ 5 \\ 1.3337 \\ 100,050 [0\%]} & \makecell{0 \\ 0 \\ 0 \\ 7,432.4 \textbf{[100\%]}} & \makecell{0.1333 \\ 0 \\ 0.3399 \\ 14,514.5 [86.7\%]} & \makecell{1.3667 \\ 0 \\ 2.7384 \\ 26,826.8 [56.7\%]} & \makecell{0 \\ 0 \\ 0 \\ 350.3 \textbf{[100\%]}} \\ \hline
  $f_{11}$ & \makecell{Avg \\ Mdn \\ SD \\ NFE [PSR]} & \makecell{3.4449e-255 \\ 2.0673e-264 \\ 0 \\ 13,413.4 \textbf{[100\%]}} & \makecell{4.3177e-240 \\ 3.7965e-242 \\ 0 \\ 3,451.7 \textbf{[100\%]}} & \makecell{7.9675e-08 \\ 1.5504e-08 \\ 1.3562e-07 \\ 21,196.2 \textbf{[100\%]}} & \makecell{2.2508e-15 \\ 2.4628e-15 \\ 1.6665e-15 \\ 21,821.8 \textbf{[100\%]}} & \makecell{3.0985e-67 \\ 1.4344e-67 \\ 6.0618e-67 \\ 5,580.6 \textbf{[100\%]}} & \makecell{0 \\ 0 \\ 0 \\ 850.8 \textbf{[100\%]}} \\ \hline
\end{tabular}

    \caption{Minimization result of benchmark functions from Table 9 \& 10 (Unimodal test functions) found by WOA, TLBO, BBO, DE, ACO and PSA, where $n = 30$ and maximum number of iterations are $1000$ }
\end{table}

\begin{table}[H]
\tiny
    \centering
\begin{tabular}{lllllllll}
\hline \hline
 &  &  FA & BA & PSO & HS & BH & GSA \\ \hline \hline
 $f_{1}$ & \makecell{Avg \\ Mdn \\ SD \\ NFE [PSR]} & \makecell{3.1585 \\ 3.2604 \\ 0.3448  \\ 2,200,356 [0\%]} & \makecell{31.4751 \\ 0.1947 \\ 49.1488 \\ 1,275,050 [40\%]} & \makecell{0.0096 \\ 2.3295e-12 \\ 0.0519 \\ 24,374.3 [96.7\%]} & \makecell{0.0158 \\ 0.0158 \\ 0.0021 50,050 [0\%]} & \makecell{2.4951 \\ 2.532 \\ 0.3481 \\ 50,050 [0\%]} & \makecell{2.8493e06 \\ 1.4669e04 \\ 83.6615e05 \\ 50,050 [0\%]} \\ \hline
$f_{2}$ & \makecell{Avg \\ Mdn \\ SD \\ NFE [PSR]} & \makecell{45.7822 \\ 41.6323 \\ 15.1427 \\ 2,170,322.5 [0\%]} & \makecell{15.8948e03 \\ 9.7789 \\ 38.9877e03 \\ 1,275,050 [0\%]} & \makecell{29.2921e04 \\ 68.9928e03 \\ 46.7798e04 \\ 50,050 [0\%]} & \makecell{10.0314 \\ 9.3557 \\ 2.9611 \\ 50,050 [0\%]} & \makecell{0.9917 \\ 0.6708 \\ 0.6884  \\ 50,050 [0\%]} & \makecell{0.6755 \\ 0.6667 \\ 0.0300 \\ 50,050 [0\%]} \\ \hline
$f_{3}$ & \makecell{Avg \\ Mdn \\ SD \\ NFE [PSR]} & \makecell{1.8942e02 \\ 1.8805e02 \\ 21.7892 \\ 2,199,033.5 [0\%]} & \makecell{6.9592e02 \\ 49.8555 \\ 12.9251e02 \\ 1,275,050 [0\%]} & \makecell{60.9412e02 \\ 19.8971 \\ 13.2010e03 \\ 50,050 [3.3\%]} & \makecell{2.7806 \\ 2.8373 \\ 0.2990 \\ 50,050 [0\%]} & \makecell{0.0054 \\ 0.0028 \\ 0.0073 \\ 50,050 [3.3\%]} & \makecell{7.272e-05 \\ 6.5751e-05 \\ 3.4463e-05 \\ 37,412.4 [76.7\%]} \\ \hline
$f_{4}$ & \makecell{Avg \\ Mdn \\ SD \\ NFE [PSR]} & \makecell{4.2587 \\ 4.314 \\ 0.5768 \\ 2,291,196.5 [0\%]} & \makecell{2.353 \\ 0.5 \\ 8.1566 \\ 887,434.8 [73.3\%]} & \makecell{27.4722 \\ 5.336 \\ 39.4775 \\ 50,050 [6.7\%]} & \makecell{\textbf{0.5000} \\ \textbf{0.5000} \\ 5.1406e-06 \\ 31,281.2 \textbf{[100\%]}} & \makecell{0.5 \\ 0.5 \\ 1.2651e-12 \\ 6,431.4 \textbf{[100\%]}} & \makecell{0.5 \\ 0.5 \\ 0 \\ 12,612.6 \textbf{\textbf{[100\%]}}} \\ \hline
$f_{5}$ & \makecell{Avg \\ Mdn \\ SD \\ NFE [PSR]} & \makecell{5.2905e02 \\ 2.5318e02 \\ 5.4441e02 \\ 2,008,748 [0\%]} & \makecell{8.9811e06 \\ 35.1230e02 \\ 20.5514e06 \\ 1,275,050 [0\%]} & \makecell{82.5440e06 \\ 35.7361e06 \\ 84.4992e06 50,050 [0\%]} & \makecell{5.6377e02 \\ 4.8958e02 \\ 3.8020e02 \\ 50,050 [0\%]} & \makecell{1.0514e02 \\ 78.0888 \\ 1.0116e02 \\ 50,050 [0\%]} & \makecell{28.8558 \\ 25.7126 \\ 13.6662 \\ 50,050 [0\%]} \\ \hline
$f_{6}$ & \makecell{Avg \\ Mdn \\ SD \\ NFE [PSR]} & \makecell{10.3257 \\ 8.7419  \\ 5.0396 \\ 1,339,726.5 [0\%]} & \makecell{94.9561e02 \\ 44.6956e02 \\ 10.4469e03 \\ 1,275,050 [0\%]} & \makecell{30.1030e03 \\ 31.4372e03 \\ 11.8775e03 \\ 50,050 [0\%]} & \makecell{61.0023e02 \\ 57.4462e02 \\ 18.8530e02 \\ 50,050 [0\%]} & \makecell{1.6389e02 \\ 1.6942e02 \\ 45.5732 \\ 50,050 [0\%]} & \makecell{1.0234e02 \\ 1.0104e02 \\ 35.8842 \\ 50,050 [0\%]} \\ \hline
$f_{7}$ & \makecell{Avg \\ Mdn \\ SD \\ NFE [PSR]} & \makecell{9.2836 \\ 7.0446  \\ 4.939 \\ 1,349,762.5 [0\%]} & \makecell{1.6889e02 \\ 54.4746 \\ 2.1129e02 \\ 1,275,050 [0\%]} & \makecell{0.0403 \\ 0.0002 \\ 0.1187 \\ 50,050 [36.7\%]} & \makecell{12.3562 \\ 12.4165 \\ 0.9441 \\ 50,050 [0\%]} & \makecell{0.0002 \\ 8.9254e-05 \\ 0.0004 \\ 50,050 [53.3\%]} & \makecell{0.4660 \\ 1.5737e-08 \\ 1.8153 \\ 31,331.3 [80\%]} \\ \hline
$f_{8}$ & \makecell{Avg \\ Mdn \\ SD \\ NFE [PSR]} & \makecell{5.9462 \\ 5.9285 \\ 0.3582 \\ 2,128,931.5 [0\%]} & \makecell{18.6205 \\ 9.8503 \\ 24.1327 \\ 1,275,050 [0\%]} & \makecell{0.0319 \\ 6.4271e-05 \\ 0.1476 \\ 36,211.2 [60\%]} & \makecell{1.188 \\ 1.1897 \\ 0.1239 \\ 50,050 [0\%]} & \makecell{2.3094e-06 \\ 4.517e-07 \\ 5.7193e-06 \\ 33,483.4 \textbf{[100\%]}} & \makecell{0.0400 \\ 1.5185e-08  \\ 0.1834 \\ 31,506.8 [93.3\%]} \\ \hline
$f_{9}$ & \makecell{Avg \\ Mdn \\ SD  \\NFE [PSR]} & \makecell{1.9569 \\ 1.9741  \\ 0.1910 \\ 1,718,742 [0\%]} & \makecell{43.1083e02 \\ 1.0844e-09 \\ 86.7370e02  \\ 314,937.3 [56.7\%]} & \makecell{2.3678e02 \\ 0.0145 \\ 7.2513e02 \\ 39,514.5 [40\%]} & \makecell{13.1538 \\ 12.9999 \\ 1.8022 \\ 50,050 [0\%]} & \makecell{1.3031e-07 \\ 7.4392e-08 \\ 1.4585e-07 \\ 32,157.1 \textbf{[100\%]}} & \makecell{1.0826e-17 \\ 9.6156e-18 \\ 4.4511e-18 \\ 16,942 \textbf{[100\%]}} \\ \hline
$f_{10}$ & \makecell{Avg \\ Mdn \\ SD \\ NFE [PSR]} & \makecell{1.3667 \\ 1 \\ 0.54671 \\ 1,731,419.5 [3.3\%]} & \makecell{2.2032e03 \\ 4.985e02 \\ 3.7426e03 \\ 1,275,050 [0\%]} & \makecell{17.1257e02 \\ 9\\ 6.0310e03 \\ 50,050 [6.7\%]} & \makecell{14.1667\\ 14\\ 1.9336\\ 50,050 [0\%]} & \makecell{0.066667\\ 0 \\ 0.24944 \\ 28,153.1 [93.3\%]} & \makecell{0 \\ 0 \\ 0 \\ 6,206.2 \textbf{[100\%]}} \\ \hline
$f_{11}$ & \makecell{Avg \\ Mdn \\ SD \\ NFE [PSR]} & \makecell{26.7122 \\ 26.2813 \\ 3.0663 \\ 2,179,059.5 [0\%]} & \makecell{7.7361e02 \\ 73.2921 \\ 12.0210e02 \\ 1,033,428 [20\%]} & \makecell{5.1003e-05 \\ 1.6173e-08 \\ 0.0002 \\ 41,241.2 [96.7\%]} & \makecell{1.777 \\ 1.8128 \\ 0.24023 \\ 50,050 [0\%]} & \makecell{2.787e-07 \\ 4.8899e-08 \\ 6.6026e-07 \\ 34,334.3 \textbf{[100\%]}} & \makecell{9.0687e-17 \\ 8.3128e-17 \\ 3.0573e-17 \\ 21,121.1 \textbf{[100\%]}} \\ \hline \hline 
\end{tabular}

    \caption{Minimization result of benchmark functions from Table 9 \& 10 (Unimodal test functions) found by FA, BA, PSO, HS, BH and GSA, where $n = 30 $ and maximum number of iterations is $1000$
}
    \label{TABLE 2 (Continued)}
\end{table}

\textbf{Multimodal functions}. As Table 1 \& 2 illustrate, PSA achieves the best values in all the multimodal test problems. In each case, there are few other algorithms whose results are similar to the results found by PSA. Exceptions are $g_{5}$, $g_{7}$ , $g_{8}$ and $g_{11}$ where the performance of PSA is strictly better than that of all other algorithms. The reported Avg, Mdn and SD results for multimodal functions have been summarized in Table 5 \& 6.

\begin{table}[H]
\tiny
    \centering
\begin{tabular}{lllllllll}
\hline \hline
 &  &   WOA & TLBO   & BBO & DE & ACO & PSA \\ \hline \hline
 $g_{1}$ & \makecell{Avg \\ Mdn \\ SD \\ NFE [PSR]} & \makecell{3.9672e-15 \\ 4.4409e-15 \\ 2.5509e-15 \\ 14,990 \textbf{[100\%]}} & \makecell{4.4409e-15 \\ 4.4409e-15 \\ 0 \\ 4,902.4 \textbf{[100\%]}} & \makecell{7.5645e-05 \\ 6.0913e-05 \\ 5.5484e-05 \\ 44,769.7 [80\%]} & \makecell{3.7728e-08 \\ 3.4195e-08 \\ 1.5908e-08 \\ 31,281.2 \textbf{[100\%]}} & \makecell{0.0310 \\ 7.9936e-15 \\ 0.1672 \\ 8,358.3 [96.7\%]} & \makecell{1.125e-15 \\ 8.8818e-16 \\ 1.2755e-15 \\ 1,226.2 \textbf{[100\%]}} \\ \hline
$g_{2}$ & \makecell{Avg\\ Mdn \\ SD \\ NFE [PSR]} & \makecell{-0.7622 \\ -0.7766 \\ 0.0559 \\ 50,050 [0\%]} & \makecell{-0.6860 \\ -0.6881 \\ 0.0338 \\ 100,050 [0\%]} & \makecell{-1 \\ -1 \\ 4.2973e-10 \\ 9,959.9 \textbf{[100\%]}} & \makecell{-0.6289 \\ -0.6105 \\ 0.0500 \\ 50,050 [0\%]} & \makecell{-0.6068 \\ -0.601 \\ 0.0275 \\ 50,050 [0\%]} & \makecell{-1  \\ -1 \\ 0  \\ 1,126.1 \textbf{[100\%]}} \\ \hline
$g_{3}$ & \makecell{Avg \\ Mdn \\ SD \\ NFE [PSR]} & \makecell{-0.9297 \\ -0.9362 \\ 0.0526 \\ 50,050 [13.3\%]} & \makecell{-0.9362 \\ -0.9362 \\ 4.9521e-17 \\ 100,050 [0\%]} & \makecell{-0.3311 \\ -0.2888 \\ 0.0838 \\ 50,050 [0\%]} & \makecell{-0.7871 \\ -0.7857 \\ 0.0058 \\ 50,050 [0\%]} & \makecell{-0.6568 \\ -0.6195 \\ 0.1099 \\ 50,050 [0\%]} & \makecell{-0.9936 \\ -1\\ 0.0191 \\ 1,651.6 [90\%]} \\ \hline
$g_{4}$ & \makecell{Avg \\ Mdn \\ SD  \\ NFE [PSR]} & \makecell{0.0020 \\ 0 \\ 0.0110 \\ 13,588.6 [96.7\%]} & \makecell{0 \\ 0 \\ 0 \\ 4,152.1 \textbf{[100\%]}} & \makecell{0.0102 \\ 0.0086 \\ 0.0100 \\ 41,016 [30\%]} & \makecell{0.0025 \\ 7.7549e-14 \\ 0.0046 \\  25,300.3 [76.7\%]} & \makecell{0.0067 \\ 0  \\ 0.0128 \\ 6,831.8 [66.7\%]} & \makecell{0 \\ 0 \\ 0 \\ 775.8 \textbf{[100\%]}} \\ \hline
$g_{5}$ & \makecell{Avg \\ Mdn \\ SD \\ NFE [PSR]} & \makecell{0.7933 \\ 0.7471 \\ 0.3219 \\ 50,050 [0\%]} & \makecell{2.091 \\ 2.1168 \\ 0.3355 \\ 100,050 [0\%]} & \makecell{0.4483 \\ 1.2967e-09 \\ 0.6559 \\ 32,532.5 [53.3\%]} & \makecell{0.1445 \\ 1.5596e-15 \\ 0.2633 \\ 35,660.6 [63.3\%]} & \makecell{1.2708 \\ 0.9534 \\ 0.9554 \\ 50,050 [0\%]} & \makecell{1.4998e-32 \\ 1.4998e-32 \\ 1.0948e-47 \\ 850.8 \textbf{[100\%]}} \\ \hline
$g_{6}$ & \makecell{Avg \\ Mdn \\ SD \\ NFE [PSR]} & \makecell{2.0075 \\ 0 \\ 3.9385 \\ 22,998 [73.3\%]} & \makecell{2.3438 \\ 3.3723 \\ 2.0551 \\ 100,050 [16.7\%]} & \makecell{6.3019 \\ 6.4489 \\ 0.8323 \\ 50,050 [0\%]} & \makecell{4.6805 \\ 4.8145 \\ 1.2038 \\ 50,050 [3.3\%]} & \makecell{12.3604  \\ 12.3783 \\ 0.2587 \\ 50,050 [0\%]} & \makecell{0 \\ 0 \\ 0 \\ 1,001 \textbf{[100\%]}} \\ \hline
$g_{7}$ & \makecell{Avg \\ Mdn \\ SD \\ NFE [PSR]} & \makecell{1.0839 \\ 1.1711 \\ 0.5197 \\ 50,050 [0\%]} & \makecell{2.6068 \\ 2.6507 \\ 0.3785 \\ 100,050 [0\%]} & \makecell{8.7806e-08 \\ 6.5568e-09 \\ 2.9285e-07 \\ 22,397.4 \textbf{[100\%]}} & \makecell{0.1208 \\ 1.526e-12 \\ 0.6456 \\ 26,952 [86.7\%]} & \makecell{0.2464 \\ 2.4591e-31 \\ 0.7697 \\ 28,528.5 [63.3\%]} & \makecell{\textbf{1.3498e-32} \\ \textbf{1.3498e-32} \\ 5.4738e-48 \\ 650.6 \textbf{[100\%]}} \\ \hline
$g_{8}$ & \makecell{Avg \\ Mdn \\ SD \\ NFE [PSR]} & \makecell{3.612e-258 \\ 9.5972e-264 \\ 0 \\ 13,363.3 \textbf{[100\%]}} & \makecell{2.4971e-239 \\ 1.0067e-240 \\ 0 \\ 4,502.2 \textbf{[100\%]}} & \makecell{3.9828e02 \\ 3.4211e02 \\ 2.9528e02 \\ 50.050 [0\%]} & \makecell{11.6561 \\ 2.7311e-07 \\ 33.188 \\ 37,837.8 [86.7\%]} & \makecell{5.1257e02 \\ 4.8337e02 \\ 4.3578e02 \\ 50,050 [16.7\%]} & \makecell{0 \\ 0 \\ 0 \\ 950.9 \textbf{[100\%]}} \\ \hline
$g_{9}$ & \makecell{Avg \\ Mdn \\ SD \\ NFE [PSR]} & \makecell{28.3982 \\ 27.9624 \\ 9.0584 \\ 50,050 [0\%]} & \makecell{66.12 \\ 66.6546 \\ 9.8457 \\ 100,050 [0\%]} & \makecell{0.0027 \\ 0.0016 \\ 0.0038 \\ 50,050 [6.7\%]} & \makecell{0.0002 \\ 3.3341e-05 \\ 0.0005 \\ 48,898.8 [70\%]} & \makecell{9.3259e-15 \\ 6.8834e-15 \\ 7.6111e-15 \\ 10,335.3 \textbf{[100\%]}} & \makecell{0 \\ 0 \\ 0 \\ 1,301.3 \textbf{[100\%]}} \\ \hline
$g_{10}$ & \makecell{Avg \\ Mdn \\ SD \\ NFE [PSR]} & \makecell{0.1199 \\ 0.0999 \\ 0.0748 \\ 50,050 [16.7\%]} & \makecell{0.0999 \\ 0.0999 \\ 0 \\ 100,050 [0\%]} & \makecell{0.8467 \\ 0.7999 \\ 0.1686 \\ 50,050 [0\%]} & \makecell{0.2247 \\ 0.2001 \\ 0.0406 \\ 50,050 [0\%]} & \makecell{0.3832 \\ 0.2999 \\ 0.2130 \\ 50,050 [0\%]} & \makecell{\textbf{0.0333} \\ \textbf{0} \\ 0.0649 \\ \textbf{26,051 [73.3\%]}} \\ \hline
$g_{11}$ & \makecell{Avg \\ Mdn \\ SD \\ NFE [PSR]} & \makecell{4.3827e-78 \\ 4.0478e-99 \\ 2.3575e-77 \\ 21,221.2 \textbf{[100\%]}} & \makecell{9.5776e-238 \\ 1.2101e-239 \\ 0 \\ 4,402.2 \textbf{[100\%]}} & \makecell{0.0010 \\ 0.0008 \\ 0.0006 \\ 50,050 [0\%]} & \makecell{2.3925e-10 \\ 1.4233e-10 \\ 3.523e-10 \\ 30,305.3 \textbf{[100\%]}} & \makecell{1.7177e-46 \\ 1.4223e-47 \\ 6.4331e-46 \\ 9,584.6 \textbf{[100\%]}} & \makecell{\textbf{-5.1284e-321} \\ \textbf{-4.6689e-321} \\ \textbf{0} \\ \textbf{1,076.1} \textbf{[100\%]}} \\ \hline
$g_{12}$ & \makecell{Avg \\ Mdn \\ SD \\ NFE [PSR]} & \makecell{-0.4961 \\ -0.0189  \\ 1.0452 \\ 50,050 [0\%]} & \makecell{-2.8297 \\ -2.8175 \\ 0.17595 \\ 100,050 [0\%]} & \makecell{-2.6667 \\ -3.5 \\ 1.1785 \\ 37,387.3 [66.7\%]} & \makecell{-0.7986 \\ -0.6578 \\ 0.3233 \\ 50,050 [0\%]} & \makecell{-0.0087 \\ -0.0059 \\ 0.0111 \\ 50,050 [0\%]} & \makecell{-3.3336 \\ -3.5 \\ 0.6225 \\ 1,526.5 [93.3\%]} \\ \hline
$g_{13}$ & \makecell{Avg \\ Mdn \\ SD \\ NFE [PSR]} & \makecell{-35.0577 \\ -35.1778 \\ 2.3959 \\ 50,050 [0\%]} & \makecell{-25.4244 \\-25.35 \\ 1.1513  \\ 100,050 [0\%]} & \makecell{-35.9776 \\ -36.1032 \\ 1.2025 \\ 50,050 [0\%]} & \makecell{-38.1295 \\ -38.2237 \\ 0.54951\\  50,050 [6.7\%]} & \makecell{-35.5849 \\ -35.632 \\ 0.9929 \\ 50,050 [0\%]} & \makecell{\textbf{-39.1662} \\ \textbf{-39.1662} \\ \textbf{0} \\ \textbf{500.5} \textbf{[100\%]}} \\ \hline
\end{tabular}
    \caption{Minimization result of benchmark functions from Table 10 \& 11 (Multimodal test functions) found by WOA, TLBO, BBO, DE, ACO and PSA, where $n =30 $and maximum number of iterations is 1000 }
    \label{TABLE 3}
\end{table}

\begin{table}[H]
\tiny
    \centering
\begin{tabular}{lllllllll}
\hline \hline
 &  &  FA & BA & PSO & HS & BH & GSA \\ \hline \hline
 $g_{1}$ & \makecell{Avg \\ Mdn \\ SD \\ NFE [PSR]} & \makecell{19.6837 \\ 19.7007 \\ 0.2346 \\ 2,289,964.5 [0\%]} & \makecell{17.8574 \\ 18.5735 \\ 1.9865 \\ 1,275,050 [0\%]} & \makecell{12.4815 \\ 13.509 \\ 5.5947 \\ 50,050 [0\%]} & \makecell{1.7832 \\ 1.8374 \\ 0.18052 \\ 50,050 [0\%]} & \makecell{0.01024 \\ 0.0018 \\ 0.0214 \\ 50,050 [13.3\%]} & \makecell{2.3768e-09 \\ 2.3432e-09 \\ 3.0133e-10 \\ 27,552.5 \textbf{[100\%]}} \\ \hline
$g_{2}$ & \makecell{Avg \\ Mdn \\ SD \\ NFE [PSR]} & \makecell{-0.6624 \\ -0.6630  \\ 0.01560  \\ 2,367,082 [0\%]} & \makecell{-0.9919 \\ -0.9998 \\ 0.0253 \\ 867,034 [53.3\%]} & \makecell{-0.8293 \\ -0.8807 \\ 0.1126 \\ 50,050 [0\%]} & \makecell{-0.9844 \\ -0.9846 \\ 0.0021 \\ 50,050 [0\%]} & \makecell{-0.9855 \\ -0.9886 \\ 0.0085 \\ 50,050 [0\%]} & \makecell{-0.9989 \\ -1 \\ 0.0060 \\ 19,744.7 [96.7\%]} \\ \hline
$g_{3}$ & \makecell{Avg \\ Mdn \\ SD \\ NFE [PSR]} & \makecell{-0.6094 \\ -0.6193 0.0352 \\ 2,235,940 [0\%]} & \makecell{-0.0536 \\ -0.0481 \\ 0.0250 \\ 1,275,050 [0\%]} & \makecell{-0.0536 \\ -0.0431 \\ 0.0321 \\ 50,050 [0\%]} & \makecell{-0.6674 \\ -0.6195 \\ 0.0739 \\ 50,050 [0\%]} & \makecell{-0.5392 \\ -0.4778 \\ 0.0702 \\ 50,050 [0\%]} & \makecell{-0.2735 \\ -0.2882 \\ 0.1233 \\ 50,050 [0\%]} \\ \hline
$g_{4}$ & \makecell{Avg \\ Mdn \\ SD \\ NFE [PSR]} & \makecell{2.2736e02 \\ 2.3281e02 \\ 39.963 \\ 1,277,054 [0\%]} & \makecell{30.6158 \\ 4.8173 \\ 58.4573 \\ 1,275,050 [16.7\%]} & \makecell{9.9293 \\ 0.2009 \\ 50.3343 \\ 50,050 [3.3\%]} & \makecell{1.1191 \\ 1.1217 \\ 0.0178 \\ 50,050 [0\%]} & \makecell{0.0180 \\ 0.0099 \\ 0.0198 \\ 50,050 [20\%]} & \makecell{2.1985 \\ 1.9595 \\ 1.0096 \\ 50,050 [0\%]} \\ \hline
$g_{5}$ & \makecell{Avg \\ Mdn \\ SD \\ NFE [PSR]} & \makecell{36.89 \\ 36.946 \\ 8.8642 \\ 2,277,246 [0\%]} & \makecell{34.0034 \\ 32.4314 \\ 10.5226 \\ 1,275,050 [0\%]} & \makecell{33.9449 \\ 23.3921 \\ 25.8672 \\ 50,050 [3.3\%]} & \makecell{0.0680 \\ 0.0696 \\ 0.0111 \\ 50,050 [0\%]} & \makecell{0.0060 \\ 5.3948e-08 \\ 0.0223 \\ 34,159.1 [90\%]} & \makecell{0.0936 \\ 4.7577e-18 \\ 0.1895 \\ 16,416.4 [73.3\%]} \\ \hline
$g_{6}$ & \makecell{Avg \\ Mdn \\ SD \\ NFE [PSR]} & \makecell{11.615 \\ 11.6374 \\ 0.3786 \\ 2,262,254 [0\%]} & \makecell{10.9257 \\ 10.9854 \\ 0.3340 \\ 1,275,050 [0\%]} & \makecell{12.1838 \\ 12.3338 \\ 0.6422 \\ 50,050 [0\%]} & \makecell{6.596 \\ 6.6573 \\ 0.4748 \\ 50,050 [0\%]} & \makecell{10.551 \\ 10.5379 \\ 0.5729 \\ 50,050 [0\%]} & \makecell{9.5538 \\ 9.2904 \\ 1.0491  \\ 50,050 [0\%]} \\ \hline
$g_{7}$ & \makecell{Avg \\ Mdn \\ SD \\ NFE [PSR]} & \makecell{0.3273 \\ 0.3317 \\0.0241 \\ 1,891,140 [0\%]} & \makecell{11.2764e06 \\ 61.1277 \\ 27.5248e06 \\ 1,275,050 [0\%]} & \makecell{51.1919e07 \\ 54.2408e07 \\ 35.0222e07 \\ 50,050 [0\%]} & \makecell{0.8414 \\ 0.7970 \\ 0.1511 \\ 50,050 [0\%]} & \makecell{0.0008 \\ 4.097e-06 \\ 0.0027 \\ 45,195.15 [83.3\%]} & \makecell{0.0015 \\ 1.1127e-18 \\ 0.0037 \\ 32,432.4 [86.7\%]} \\ \hline
$g_{8}$ & \makecell{Avg \\ Mdn \\ SD \\ NFE [PSR]} & \makecell{2.1681e03  \\ 20.2030e02 \\ 6.2277e02 \\ 2,002,219 [0\%]} & \makecell{6.4985e03 \\ 64.0538e02 \\ 10.4159e02 \\ 1,275,050 [0\%]} & \makecell{4.2626e03 \\ 39.7341e02 \\ 17.5046e02 \\ 50,050 [0\%]} & \makecell{1.5742e02 \\ 1.4514e02 \\ 42.646 \\ 50,050 [0\%]} & \makecell{1.8685e03 \\ 19.3556e02 \\ 10.0055e02 \\ 50,050 [0\%]} & \makecell{2.3091e02 \\ 2.0661e02 \\ 1.5733e02 \\ 50,050 [0\%]} \\ \hline
$g_{9}$ & \makecell{Avg \\ Mdn \\ SD \\ NFE [PSR]} & \makecell{56.3622 \\ 55.6749 \\ 10.7762 \\ 2,179,551 [0\%]} & \makecell{29.5766e02 \\ 18.3567 \\ 97.9610e02 \\ 1,275,050 [0\%]} & \makecell{85.7735e03 \\ 68.7687e03 \\ 77.8947e03 \\ 50,050 [0\%]} & \makecell{14.8272 \\ 14.6629 \\ 1.8852 \\ 50,050 [0\%]} & \makecell{0.0072 \\ 0.0005 \\ 0.0156 \\ 49,899.8 [16.7\%]} & \makecell{4.4172 \\ 2.1003\\ 4.9525 \\ 50,050 [33.3\%]} \\ \hline
$g_{10}$ & \makecell{Avg \\ Mdn \\ SD \\ NFE [PSR]} & \makecell{0.3873 \\ 0.2999 \\ 0.1498 \\ 1,293,717.5 [0\%]} & \makecell{7.3632 \\ 6.6499 \\ 4.7592 \\ 1,275,050 [0\%]} & \makecell{7.33 \\ 6.3353 \\ 2.8086 \\ 50,050 [0\%]} & \makecell{1.2166 \\ 1.2168 \\ 0.1583 \\ 50,050 [0\%]} & \makecell{0.6899 \\ 0.6999 \\ 0.1012 \\ 50,050 [0\%]} & \makecell{1.586 \\ 1.6893 \\ 0.3420 \\ 50,050 [0\%]} \\ \hline
$g_{11}$ & \makecell{Avg \\ Mdn \\ SD \\ NFE [PSR]} & \makecell{41.7632 \\ 41.6016 \\ 4.5022 \\ 1,651,058 [0\%]} & \makecell{59.2780e04 \\ 58.1649e04 \\ 11.3170e04  \\ 1,275,050 [0\%]} & \makecell{18.1836e02 \\ 2.1673e02 \\ 51.0800e02 \\ 50,050 [0\%]} & \makecell{5.5315e02 \\ 5.2091e02 \\ 1.2103e02 \\ 50,050 [0\%]} & \makecell{0.0008 \\ 0.0004 \\ 0.0009 \\ 48,823.8 [3.3\%]} & \makecell{2.4097e-16 \\ 2.4479e-16 \\ 6.1442e-17 \\ 21,221.2 \textbf{[100\%]}} \\ \hline
$g_{12}$ & \makecell{Avg \\ Mdn \\ SD \\ NFE [PSR]} & \makecell{-0.9480 \\ -0.9452 \\ 0.0835 \\ 2,277,505.5 [0\%]} & \makecell{-3.2293 \\ -3.2937 \\ 0.2281 \\ 1,275,050 [0\%]} & \makecell{-0.0077 \\ -0.0060 \\ 0.0086 \\ 50,050 [0\%]} & \makecell{-0.4392 \\ -0.4262 \\ 0.0974 \\ 50,050 [0\%]} & \makecell{-2.0039 \\ -2.4114 \\ 1.1432 \\ 50,050 [0\%]} & \makecell{-3.0833 \\ -3.5 \\ 0.9317 \\ 34,434.4 [83.3\%]} \\ \hline
$g_{13}$ & \makecell{Avg \\ Mdn \\ SD \\ NFE [PSR]} & \makecell{-33.0316 \\ -32.9629 \\ 0.8953 \\ 2,226,590 [0\%]} & \makecell{-33.6979 \\ -34.9251 \\ 3.1559  \\ 1,275,050 [0\%]} & \makecell{-18.7046 \\ -17.2052 \\ 5.4533 \\ 50,050 [0\%]} & \makecell{-39.1494 \\ -39.149 \\ 0.0029 \\ 50,050 [0\%]} & \makecell{-33.7709\\ -33.9049 \\ 1.2115 \\ 50,050 [0\%]} & \makecell{-33.9827 \\ -33.9827 \\ 1.0466 \\ 50,050 [0\%]} \\ \hline
\end{tabular}
    \caption{Minimization result of benchmark functions from Table 10 (Multimodal test functions) found by FA, BA, PSO, HS, BH and GSA, where $n =30 $and maximum number of iterations is 1000 }
    \label{TABLE 3 continue}
\end{table}

By comparing these obtained results with the best-so-far values of average, median and standard deviation of PSA, we can see that PSA provides the best results except for the standard deviations related to $g_{1}$ , $g_{3}$ and $g_{10}$ (compared with TLBO) and $g_{12}$ (compared with PSO). For other results, there are few other algorithms (in comparison with the unimodal functions) with the same performance as the PSA. Also, for $g_{5}$ , $g_{7}$ , $g_{9}$ and $g_{13}$ , PSA provides strictly better results than all the algorithms in all the criteria Avg, Mdn and SD.\\
Table 3 also include the number of function evaluation (NFE) and the percentage of successful runs (PSR) for the algorithms. Based on these tables, PSA has the maximum successful rates and uses the minimum number of evaluation to reach the required accuracy. Therefore, the performance of PSA is shown to be relatively good, as it achieves the best PSR (=100) in 10 out of 13 test problems, $PSR=90\%$ for $g_{3}$ , $PSR=73.3\%$ for $g_{10}$ and $PSR=93.3\%$ for $g_{12}$.\\
These results prove the good convergence rate, accuracy and stability of PSA in handling the multimodal functions.\\

\textbf{Shifted functions}. As mentioned previously, functions $\widetilde{h_{i}}$ have resulted from basic test function $h_{i}$, listed in table A.5 , A.6, by shifting their global optimums. More precisely, $\widetilde{h_{i}}(x) = h_{i}(x-\textbf{1}), i = 1,2,...,11$, where $\textbf{1}$ is an n-dimensional vector with each component equal to one. Therefore, the optimal values will remain unchanged. In other words, if $x^{*}$ denotes the optimum of $\widetilde{h_{i}}$ then $x^{*} - \textbf{1}$ is the optimum of $h_{1}$, and $\widetilde{h_{i}}(x^{}) = h_{i}(x^{} - \boldsymbol{1})$. Table 7 and 8 present the basic benchmark function $\widetilde{h_{i}}, i=1,2,...,11$. Some of the test functions are unimodal and others are multimodal functions. Comparison of the results of best values found by PSA and those obtained using other algorithms show that PSA is the only best performing algorithm for all the shifted test problems; results of which are presented in Table 1 and 2. In this respect only for $\widetilde{h_{3}}$, ACO, and GSA perform similarly well.\\

On the shifted test functions (especially, multimodal ones with a large number of local optima), finding good solutions and escaping from local optima is very hard. However, by comparing the indices Avg, Mdn and SD which are summarized in Table 4, it can be seen that PSA exhibits significant performance, and provides much better results in comparison to other algorithms on these functions. Although, the only exceptions are Avg and SD for $\widetilde{h_{1}}$ obtained by ACO, Avg and SD for $\widetilde{h_{3}}$ obtained by BBO, SD for $\widetilde{h_{1}}$ obtained by TLBO and SD for $\widetilde{h_{10}}$ obtained by FA. \\
Furthermore, PSA achieved the best PSR (=100) for 8 out of 11 test problems. Also, it attained PSR=96.7\% for $\widetilde{h_{3}}$, PSR=53.3\% for $\widetilde{h_{1}}$ and PSR = 63.3\% for $\widetilde{h_{10}}$. Additionally, for $\widetilde{h_{1}}$, $\widetilde{h_{9}}$, and $\widetilde{h_{10}}$ , PSA is the only algorithm achieving the required accuracy. All other algorithms in this category require many more function evaluations in order to reach the required accuracy. Finally, like the unimodal functions ( $f_{i}$ ) and multimodal functions ( $g_{i}$ ), PSA used the minimum number of function evaluation for solving shifted functions $\widetilde{h_{i}}$ in all the test problems. To sum up, PSA is capable of handling these types of problems very effectively.

\begin{table}[H]
\tiny
    \centering
\begin{tabular}{lllllllll}
\hline \hline
 &  &   WOA & TLBO   & BBO & DE & ACO & PSA \\ \hline \hline
 $\widetilde{h_{1}}$ & \makecell{Avg \\ Mdn \\ SD \\ NFE [PSR]} & \makecell{5.2305\\ 4.7087 \\ 2.8235 \\ 50,050 [0\%]} & \makecell{8.8366 \\ 8.4279 \\ 1.5058 \\ 100,050 [0\%]} & \makecell{0.0014 \\ 0.0002 \\ 0.0032 \\ 50,050 [40\%]} & \makecell{0.0012 \\ 4.3305e-05 \\ 0.0027 \\ 44,869.8 [63.3\%]} & \makecell{2.3997e-15 \\ 1.4764e-15 \\ 2.3495e-15 \\ 7,157.1 \textbf{[100\%]}} & \makecell{5.3308e-15 \\ 1.5613e-16 \\ 1.4228e-14 \\ 925.9 \textbf{[100\%]}} \\ \hline
$\widetilde{h_{2}}$ & \makecell{Avg \\ Mdn \\ SD \\ NFE [PSR]} & \makecell{32.9372e05 \\ 31.5471e05 \\ 20.1215e05 \\ 50,050 [0\%]} & \makecell{17.9546e06 \\ 18.3487e06 \\ 29.8458e05 \\ 100,050 [0\%]} & \makecell{0.1605 \\ 0.0625 \\ 0.3014 \\ 50,050 [0\%]} & \makecell{3.2173e-08 \\ 1.2392e-08 \\ 4.4412e-08 \\ 38,763.7 \textbf{[100\%]}} & \makecell{3.858e-25 \\ 3.0199e-25 \\ 3.2807e-25 \\ 10,635.6 \textbf{[100\%]}} & \makecell{0 \\ 0 \\ 0 \\ 1,051 \textbf{[100\%]}} \\ \hline
$\widetilde{h_{3}}$ & \makecell{Avg \\ Mdn  \\ SD \\ NFE [PSR]} & \makecell{31.3882 \\ 31.1355 \\ 6.0552 \\ 50,050 [0\%]} & \makecell{73.6256 \\ 74.6363 \\ 10.1223 \\ 100,050 [0\%]} & \makecell{2.6121e-06 \\ 6.5957e-07 \\ 4.3797e-06 \\ 30,205.2 \textbf{[100\%]}} & \makecell{0.6924 \\ 1.0498 \\ 0.7700 \\ 38,288.2 [46.7\%]} & \makecell{2.2209 \\ 2.0996 \\ 1.4975 \\ 28,453.4 [10\%]} & \makecell{1.0092 \\ 0 \\ 5.4348 \\ 900.9 [96.7\%]} \\ \hline
$\widetilde{h_{4}}$ & \makecell{Avg\\ Mdn\\ SD\\ NFE [PSR]} & \makecell{25.2863\\ 10.6199\\ 37.5031\\ 50,050 [0\%]} & \makecell{3.8960e02\\ 3.7391e02\\ 1.2054e02\\ 100,050 [0\%]} & \makecell{1.0074e-12\\ 3.0819e-15\\ 4.4873e-12\\ 14,764.7 \textbf{[100\%]}} & \makecell{1.5526e-27\\ 2.092e-28 \\ 3.7287e-27\\ 17,792.8 \textbf{[100\%]}} & \makecell{1.3837e-60\\ 8.0978e-62\\ 5.4722e-60\\ 5,130.1 \textbf{[100\%]}} & \makecell{0\\ 0\\ 0\\ 875.9 \textbf{[100\%]}} \\ \hline
$\widetilde{h_{5}}$ & \makecell{Avg \\ Mdn \\ SD \\ NFE [PSR]} & \makecell{0.0030\\ 0.0003 \\ 0.0105 \\ 50,050 [26.7\%]} & \makecell{19.0203 \\ 19.1642 \\ 4.8144 \\100,050 [0\%]} & \makecell{1.0772e-17 \\ 1.3193e-18\\ 2.1481e-17\\ 23,873.8 \textbf{[100\%]}} & \makecell{1.349e-06\\ 2.2837e-07\\ 3.9214e-06\\ 22,597.6 \textbf{[100\%]}} & \makecell{6.5441e-08\\ 3.2606e-09\\ 1.3092e-07\\ 2,602.6 \textbf{[100\%]}} & \makecell{3.0601e-76\\ 1.6005e-94\\ 1.6479e-75\\ 250.2 \textbf{[100\%]}} \\ \hline
$\widetilde{h_{6}}$ & \makecell{Avg \\ Mdn \\ SD \\ NFE [PSR]} & \makecell{1.2360e03 \\ 7.9364e02 \\ 1.0787e03\\ 50,050 [0\%]} & \makecell{7.7749e03\\ 6.9387e03\\ 3.1597e03\\ 100,050 [0\%]} & \makecell{1.6365\\ 0.2409 \\3.5795 \\50,050 [0\%]} & \makecell{2.887e-12\\ 1.9148e-12\\ 3.0327e-12\\ 30,530.5 \textbf{[100\%]}} & \makecell{4.7861e-28\\ 2.3802e-28\\ 5.4757e-28\\ 7,707.7 \textbf{[100\%]}} & \makecell{0\\ 0\\ 0\\ 1,251.2 \textbf{[100\%]}} \\ \hline
$\widetilde{h_{7}}$ & \makecell{Avg\\ Mdn\\ SD \\NFE [PSR]} & \makecell{28.9636\\ 28.94\\ 0.7954 \\50,050 [0\%]} & \makecell{29.9497\\ 29.9506\\ 0.0205\\ 100,050 [0\%]} & \makecell{31.3412\\ 31.8387 \\ 9.4835\\ 50,050 [0\%]} & \makecell{16.0096e01\\ 16.4263e01 \\ 21.8432 \\50,050 [0\%]} & \makecell{34.3924\\ 32.3361\\ 10.2233\\ 50,050 [0\%]} & \makecell{19.8992\\ 0 \\ 30.1786 \\ 50,050 [53.3\%]} \\ \hline
$\widetilde{h_{8}}$ & \makecell{Avg \\ Mdn \\ SD \\ NFE [PSR]} & \makecell{0.0887\\ 0.0003\\ 0.3732\\ 50,050 [13.3\%]} & \makecell{16.4095\\ 17.0282\\ 4.9455\\ 100,050 [0\%]} & \makecell{5.8166e-09\\ 3.2843e-10\\ 1.1784e-08\\ 19,694.7 \textbf{[100\%]}} & \makecell{2.7032e-06\\ 3.1885e-15\\ 1.2511e-05\\ 6,031 \textbf{[100\%]}} & \makecell{8.2173e-34\\ 7.5965e-65 \\ 3.0746e-33 \\ 2,577.6 \textbf{[100\%]}} & \makecell{4.1087e-34 \\ 0 \\ 2.2126e-33 \\ 325.3 \textbf{[100\%]}} \\ \hline
$\widetilde{h_{9}}$ & \makecell{Avg \\Mdn \\ SD  \\NFE [PSR]} & \makecell{78.4202\\ 63.2084 \\ 56.8427 \\ 50,050 [0\%]} & \makecell{82.3708 \\ 71.655 \\ 38.5358 \\ 100,050 [0\%]} & \makecell{13.1452 \\ 11.8146 \\ 10.3958 \\ 50,050 [0\%]} & \makecell{2.0126 \\ 1.7747 \\ 2.3634 \\ 50,050 [0\%]} & \makecell{4.0905 \\ 2.0264\\ 4.1748 \\ 50,050 [0\%]} & \makecell{\textbf{0}\\ \textbf{0}\\ 0\\ 650.6 \textbf{[100\%]}} \\ \hline
$\widetilde{h_{10}}$ & \makecell{Avg \\ Mdn \\ SD \\ NFE [PSR]} & \makecell{0.4668\\ 0.4790\\ 0.0923\\ 50,050 [0\%]} & \makecell{0.5832\\ 0.5808\\ 0.0361\\ 100,050 [0\%]} & \makecell{0.3545\\ 0.3475\\ 0.0617\\ 50,050 [0\%]} & \makecell{0.6079\\ 0.6187\\ 0.0528\\ 50,050 [0\%]} & \makecell{0.4419\\ 0.5494\\ 0.2511 \\50,050 [0\%]} & \makecell{\textbf{0.0772} \\ \textbf{0} \\ \textbf{0.1198} \\ \textbf{25,175.1 [63.3\%]}} \\ \hline
$\widetilde{h_{11}}$ & \makecell{Avg\\ Mdn \\ SD\\ NFE [PSR]} & \makecell{4.7887e02 \\ 4.5577e02  \\ 81.2657 \\ 50,050 [0\%]} & \makecell{42.2103 \\ 35.6207 \\ 17.8737 \\ 100,050 [0\%]} & \makecell{31.0847 \\ 28.4182 \\ 13.8752 \\ 50,050 [0\%]} & \makecell{1.0199 \\ 0.7077 \\ 0.7690 \\ 50,050 [0\%]} & \makecell{0.0023 \\ 0.0001 \\ 0.0066 \\ 50,050 [33.3\%]} & \makecell{\textbf{0} \\ \textbf{0}\\ \textbf{0}\\ \textbf{3,753.7 [100\%]}} \\ \hline
\end{tabular}
    \caption{Minimization result of benchmark functions in Table 11 (shifted functions) found by WOA, TLBO, BBO, DE, ACO and PSA, where $n=30$ and maximum number of iterations is 1000}
    \label{TABLE 4}
\end{table}

\begin{table}[H]
\tiny
    \centering
\begin{tabular}{lllllllll}
\hline \hline
 &  &  FA & BA & PSO & HS & BH & GSA \\ \hline \hline
 $\widetilde{h_{1}}$ & \makecell{Avg \\ Mdn \\ SD \\ NFE [PSR]} & \makecell{16.5535\\ 16.7776 \\ 2.2543 \\ 2,237,292 [0\%]} & \makecell{7.4923 \\ 7.3576 \\ 3.8495 \\ 1,275,050 [0\%]} & \makecell{5.8091 \\ 0.9473 \\ 14.0008 \\ 50,050 [0\%]} & \makecell{0.0004 \\ 1.3989e-09 \\ 0.0010 \\ 26,301.3 [80\%]} & \makecell{0.0626 \\ 0.0269 \\ 0.1056 \\ 50,050 [0\%]} & \makecell{0.0004 \\ 1.3989e-09 \\ 0.0010 \\ 26,301.3 [80\%]} \\ \hline
$\widetilde{h_{2}}$ & \makecell{Avg \\ Mdn \\ SD \\ NFE [PSR]} & \makecell{18.8133e05 \\ 19.0087e05 \\ 1.5070e05 \\ 1,714,698 [0\%]} & \makecell{36.3748e08 \\ 86.7148e02 \\ 73.7024e08 \\ 1,275,050 [0\%]} & \makecell{13.1364e08 \\ 2.9052 \\ 66.8277e08 \\ 50,050 [3.3\%]} & \makecell{98.8708 \\ 44.1272 \\ 107.579 \\ 50,050 [0\%]} & \makecell{4.2282e02 \\ 1.9881e02 \\ 5.7525e02 \\ 50,050 [0\%]} & \makecell{98.8708 \\ 44.1272 \\ 1.0758e02 \\ 50,050 [0\%]} \\ \hline
$\widetilde{h_{3}}$ & \makecell{Avg \\ Mdn \\ SD \\ NFE [PSR]} & \makecell{23.0051 \\ 23.2354 \\ 1.3196 \\ 1,964,025 [0\%]} & \makecell{29.39376e02 \\ 1.0472e02 \\ 62.6571e02 \\ 1,275,050 [0\%]} & \makecell{3.7845e02 \\ 4.877 \\ 18.2722e02 \\ 50,050 [3.3\%]} & \makecell{0.61947 \\ 0.4414 \\ 0.6424 \\ 35,985.9 [40\%]} & \makecell{0.3648 \\ 0.0040 \\ 0.4921 \\ 45,570.5 [30\%]} & \makecell{0.6195 \\ 0.4414 \\ 0.6424 \\ 35,985.9 [40\%]} \\ \hline
$\widetilde{h_{4}}$ & \makecell{Avg \\ Mdn \\ SD \\ NFE [PSR]} & \makecell{4.2601 \\ 4.2018 \\ 0.7055 \\ 1,668,421 [0\%]} & \makecell{13.5384e06 \\ 79.5771 \\ 49.7340e06 \\ 1,275,050 [43.3\%]} & \makecell{16.2622e04 \\ 1.4027e-07 \\ 70.6541e04 \\ 39,939.9 [60\%]} & \makecell{6.9258e-35 \\ 5.593e-35 \\ 4.0983e-35 \\ 11,486.5 [\textbf{100\%}]} & \makecell{4.249e-11 \\ 8.1923e-13 \\ 1.922e-10 \\ 25,850.8 [\textbf{100\%}]} & \makecell{6.9258e-35 \\ 5.593e-35 \\ 4.0983e-35 \\ 11,486.5 [\textbf{100\%}]} \\ \hline
$\widetilde{h_{5}}$ & \makecell{Avg \\ Mdn \\ SD \\ NFE [PSR]} & \makecell{0.1107 \\ 0.1099 \\ 0.0273 \\ 2,266,519 [0\%]} & \makecell{3.66 \\ 0.6361 \\ 6.4116 \\ 1,275,050 [20\%]} & \makecell{1.3524 \\ 0.0442 \\ 3.1091 \\ 50,050 [16.7\%]} & \makecell{1.5744e-52 \\ 3.3009e-53 \\ 4.4847e-52 \\ 9,959.9 [\textbf{100\%}]} & \makecell{0.3703 \\ 0.3769 \\ 0.0355 \\ 50,050 [0\%]} & \makecell{1.5744e-52 \\ 3.3009e-53 \\ 4.4847e-52 \\ 9,959.9 [\textbf{100\%}]} \\ \hline
$\widetilde{h_{6}}$ & \makecell{Avg \\ Mdn \\ SD \\ NFE [PSR]} & \makecell{4.2740e04 \\ 3.8753e04 \\ 1.4213e04 \\ 1,288,368 [0\%]} & \makecell{4.0173e05 \\ 4.5656e04 \\ 1.2319e06 \\ 1,275,050 [0\%]} & \makecell{1.1783e-07 \\ 6.7044e-10 \\ 4.6328e-07 \\ 33,708.7 [\textbf{100\%}]} & \makecell{9.6724e02 \\ 9.0891e02 \\ 5.4377e02 \\ 50,050 [0\%]} & \makecell{4.6072e02 \\ 3.3705e02 \\ 4.1809e02 \\ 50,050 [0\%]} & \makecell{9.6724e02 \\ 9.0891e02 \\ 5.4377e02 \\ 50,050 [0\%]} \\ \hline
$\widetilde{h_{7}}$ & \makecell{Avg \\ Mdn \\ SD \\ NFE [PSR]} & \makecell{1.9142e02 \\ 1.9346e02 \\ 19.0159 \\ 2,281,812.5 [0\%]} & \makecell{1.6481e02 \\ 1.5773e02 \\ 41.6145 \\ 1,275,050 [0\%]} & \makecell{1.6663e02 \\ 1.4775e02 \\ 1.0037e02 \\ 50,050 [0\%]} & \makecell{19.9655 \\ 18.9042 \\ 5.744 \\ 50,050 [0\%]} & \makecell{51.9033 \\ 51.8105 \\ 11.1337 \\ 50,050 [0\%]} & \makecell{19.9655 \\ 18.9042 \\ 5.744 \\ 50,050 [0\%]} \\ \hline
$\widetilde{h_{8}}$ & \makecell{Avg \\ Mdn \\ SD \\ NFE [PSR]} & \makecell{0.0075\\ 0.0068\\ 0.0035 \\ 2,221,983.5 [0\%]} & \makecell{9.1415 \\ 0.0425 \\ 21.9564 \\ 1,275,050 [36.7\%]} & \makecell{0.1180 \\ 1.0283e-06\\ 0.5612 \\ 43,018 [76.7\%]} & \makecell{2.2962e-11 \\ 5.1698e-12 \\ 6.0402e-11 \\ 12,862.8 \textbf{[100\%]}} & \makecell{0.0713 \\ 0.0727 \\ 0.0096 \\ 50,050 [0\%]} & \makecell{2.2962e-11\\ 5.1698e-12\\ 6.0402e-11 \\ 12,862.8 \textbf{[100\%]}} \\ \hline
$\widetilde{h_{9}}$ & \makecell{Avg\\ Mdn\\ SD\\ NFE [PSR]} & \makecell{2.1263e02\\ 2.1927e02 \\ 27.6932 \\ 2,293,639.5 [0\%]} & \makecell{1.0686e02\\ 41.6048 \\1.3891e02 \\ 1,275,050 [0\%]} & \makecell{98.1528 \\ 17.3008\\ 1.9507e02\\ 50,050 [0\%]} & \makecell{14.1375 \\ 13.2904 \\ 12.947\\ 50,050 [0\%]} & \makecell{11.6503 \\ 12.8832 \\ 7.4511 \\ 50,050 [0\%]} & \makecell{14.1375 \\ 13.2904 \\ 12.947 \\ 50,050 [0\%]} \\ \hline
$\widetilde{h_{10}}$ & \makecell{Avg \\ Mdn \\ SD \\ NFE [PSR]} & \makecell{0.6116 \\ 0.6115 \\ \textbf{0.0288} \\ 2,335,961.5 [0\%]} & \makecell{0.4698 \\ 0.4690 \\ 0.0427 \\ 1,275,050 [0\%]} & \makecell{0.3232 \\ 0.301 \\ 0.1058 \\ 50,050 [0\%]} & \makecell{0.7545 \\ 0.7691 \\ 0.0541 \\ 50,050 [0\%]} & \makecell{0.3128 \\ 0.3183 \\ 0.0378 \\ 50,050 [0\%]} & \makecell{0.7545 \\ 0.7691\\ 0.0541\\ 50,050 [0\%]} \\ \hline
$\widetilde{h_{11}}$ & \makecell{Avg \\ Mdn \\ SD \\ NFE [PSR]} & \makecell{3.0594\\ 3.0853\\ 0.3933\\ 2,013,768 [0\%]} & \makecell{12.148 \\ 0.0143\\ 35.8589 \\ 1,275,050 [36.7\%]} & \makecell{8.15971e02 \\ 7.4653e02  \\ 3.8495e02  \\ 50,050 [0\%]} & \makecell{97.1866 \\ 94.3317 \\ 24.1285 \\ 50,050 [0\%]} & \makecell{0.2276 \\ 0.1820 \\ 0.1451 \\ 50,050 [0\%]} & \makecell{97.1866\\ 94.3317 \\ 24.1285 \\ 50,050 [0\%]} \\ \hline
\end{tabular}
    \caption{MINIMIZATION RESULT OF BENCHMARK FUNCTIONS IN TABLE 11 (SHIFTED FUNCTIONS) FOUND BY FA, BA, PSO, HS, BH AND GSA, WHERE $n =30 $ AND MAXIMUM NUMBER OF ITERATIONS is 1000}
    \label{TABLE 4 continue }

\end{table}

\section{Conclusion}
In this study, a novel population-based meta-heuristic optimization algorithm called Perfectionism Search Algorithm (PSA), was introduced. PSA was designed based on simulation of the most popular model and psychological aspects of perfectionism, initially proposed by Hewitt and Flett. This research aimed at presenting the idea, theoretical and applied aspects of the algorithm to continuous domains, as well as an implementation that performs well on standard benchmark test problems. In order to evaluate the efficiency of the introduced algorithm, it was applied to a set of various standard benchmark functions by solving 35 mathematical nonlinear optimization problems. These test functions were categorized into three classes of unimodal, multimodal, and shifted functions (which resulted from some unimodal and multimodal basic functions). These benchmark functions were used to make two kinds of comparisons, based on different termination criteria: i) the maximum iteration and ii) the number of function evaluations needed to reach an required accuracy. Also, the performance of PSA was compared against a substantial number of famous algorithms and approaches such as: ACO, PSO, TLBO, BBO, DE, HS, FA, WOA, BA, GSA and BH. The results obtained by PSA provided superior results in global optima achievement in comparison to the other approaches. Also, from the best-so-far results of average, median and standard deviation obtained from the execution of PSA, it was observed that this algorithm has a high stability and capability of finding high quality solutions with fast convergence rate. Moreover, PSA has high percentage of successful runs in satisfying the accuracy requirements and uses the minimum number of function evaluation in all the test problems. The optimization results demonstrated that PSA is an efficient algorithm to deal with complex optimization problems and may, therefore, be considered as a competitive approach when compared to other renowned meta- heuristic algorithms. Finally, further research into the performance of PSA on the constrained optimization problems will be the subject of future publications.

\section{Appendix A}

Table 9 (unimodal functions), Table 10 (multimodal functions) and Table 11 (shifted functions) where n denotes the dimension of a given function.

\begin{table}[H]
\footnotesize
    \centering
    \begin{tabular}{lll}
    \hline \hline
         Name & Function & Interval \\ \hline
        Brown & $f_{1}(x) = \sum_{i=1}^{n-1}((x_{i}^{2})^{(x_{i+1}^{2} + 1)} + (x_{i+1}^{2})^{(x_{i}^{2} + 1)})$ & $[-1,4]^{n}$ \\ \hline
        Dixon-Price & $f_{2}(x) = (x_{1} - \textbf{1})^{2} + \sum_{i=2}^{n}i(2x_{i}^{2}-x_{i-1})^{2}$ & $[-10,10]^{n}$ \\ \hline
        Powell singular 2 & $f_{3}(x) = \sum_{i=2}^{n-2}(x_{i-1} + 10x_{i})^{2} + 5(x_{i+1} - x_{i+2})^{2} + (x_{i} - 2x_{i+1})^{4} + 10(x_{i-1} x_{i+2})^{4}$ & $[-4,5]^{n}$\\ \hline
        Quartic & $f_{4}(x) = \sum_{i=1}^{n}ix^{4} + Random([0,))$ & $[-1.28, 1.28]^{n}$\\ \hline
        Rosenbrock & $f_{5}(x) = \sum_{i=1}^{n-1}(100(x_{i+1} - x_{i}^{2}))^{2}+(x_{i}-1)^{2}$ & $[-30,30]^{n}$ \\ \hline
        Schwefel 1.2 & $f_{6}(x) = \sum_{i=1}^{n} (\sum_{j=1}^{i}x_{j})^{2}$ & $[-100,100]^{n}$ \\ \hline
        Schwefel 2.20 & $f_{7}(x) = \sum_{i=1}^{n} \vert x_{i} \vert$ & $[-100,100]^{n}$ \\ \hline
        Schwefel 2.22 & $f_{8}(x) = \sum_{i=1}^{n} \vert x_{i} \vert + \Pi_{i=1}^{n}\vert x_{i}$ & $[-10,10]^{n}$ \\ \hline
        Sphere & $f_{9}(x) = \sum_{i=1}^{n} x_{i}^{2}$ & $[-100,100]^{n}$ \\ \hline
        Step 2 & $f_{10}(x) = \sum_{i=1}^{n} ([x_{i}+0.5])^{2}$ & $[-10,10]^{n}$ \\ \hline
        Sum squares & $f_{11}(x) = \sum_{i=1}^{n} ix_{i}^{2} $ & $[-10,10]^{n}$ \\ \hline \hline
    \end{tabular}
    \caption{UNIMODAL TEST FUNCTIONS USED IN THE EXPERIMENTAL STUDY}
\end{table}

In table 9, most of the test functions being minimized have a global minimum at $x^{*} = \textbf{0}$ (where $\textbf{0}$ is an n-dimensional zero vector) with optimal value equal zero (i.e., $f_{x^{*}} = 0$). Exceptions are as follows: Dixon-Price with $f_{2}(x^{*}) = 0$ at $[2^{(2^{1-i}-1)}]_{i=1}^{n}$; Quartic with $f_{2}(x^{*}) = Random[0.1)$ (in this paper we set $Random[0.1) = 0.5$) at $x^{*} = \textbf{0}$; Rosenbrock with $f_{5}(x^{*})=0$ at $x^{*} = [1]^{n}$; Step 2 with $f_{10}x^{x^{*}}$ at $x_{i}^{*} \in [-0.5,0.5) (i = 1,2,...,n)$. \\
In table 10, forfunctions $g_{1}$ , $g_{4}$ , $g_{6}$ , $g_{8}$ , $g_{10}$ and $g_{11}$, the global minimums are located at $x^{*} = \textbf{0}$ with $f(x^{*}) = 0$. Exceptions are as follows: Deb 1 with $g_{2}(x^{*}) = -1$ at $x_{i}^{*} = \frac{(2t_{i}+1)}{2}$ where $t \in (1,2,...,n)$; Drop wave with $g_{3}(x^{*}) = -1$ at $x^{*} = [0]^{n}$; Levy with $g_{5}(x^{*}) = 0$ at $x^{*} = [1]^{n}$; Penalty 2 with $g_{7}(x^{*}) = 0$ at $x^{*} = [1]^{n}$; Quintic with $g_{9}(x^{*}) = 0$ at $x^{*} = [-1]^{n}$ or $x^{*} = [2]^{n}$; Sinusoidal with $g_{12}(x^{*}) = -3.5$ at $x^{*} = [180]^{n}$ and and Styblinski–Tang with $g_{13}(x^{*}) = -39.16599$ at $x^{*} = [-2.903534]^{n}$. \\
In Table 11, all the test functions have a global minimum at $x^{*} = \textbf{0}$ with $f(x^{*}) = 0$. 

\begin{table}[H]
    \centering
    \begin{tabular}{lll}
    \hline \hline
    Name & Function & Initial Interval \\ \hline
    Ackley & $g_{1}(x) = -20 \exp(-0.20 \sqrt{\frac{1}{n}\sum_{i=1}^{n}x_{i}^{2}}) - \exp(\frac{1}{n}\sum_{i=1}^{n}\cos(2 \pi x_{i})) + 20 + \epsilon$ & $[-32,32]^{n}$ \\ \hline
    Deb 1 & $g_{2}(x) = -\frac{1}{n}\sum_{i=1}^{n}\sin^{6}(5\pi x_{i})$ & $[-1,1]^{n}$ \\ \hline
    Drop wave & $g_{3}(x) = - \frac{1+ \cos(12 \sqrt{\sum_{i=1}^{n}x_{i}^{2}})}{2+0.5\sum_{i=1}^{n}x_{i}^{2}}$ & $[-5.12,5.12]^{n}$ \\ \hline
    Griewank & $g_{4}(x) = 1 + (\frac{1}{4000})\sum_{i=1}^{n}x_{i}^{2} - \Pi_{i=1}^{n}\cos(\frac{x_{i}}{\sqrt{i}})$ & $[-100,100]^{n}$ \\ \hline
    Levy & \makecell{$g_{5}(x) = \sin^{2}(\pi w_{1}) + \sum_{i=1}^{n}(w_{i}-1)^{2}[1 + 10 \sin^{2}(\pi w_{i} + 1)] + (w_{n} - 1)^{2}$ \\ $w_{i} = 1 + \frac{x_{i} -1 }{4}$} & $[-10,10]^{n}$ \\ \hline
    Pathological & $g_{6}(x) = \sum_{i=1}^{n}(0.5 + \frac{\sin^{2}(\sqrt{100 x_{i}^2 + x_{i+1}^{2}}) - 0.5}{1 + 0.001(x_{i}-x_{i+1})^{4}})$ \\ \hline
    Penalty2 & \makecell{$g_{7}(x) = 0.1\{ \sin^{2}(3 \pi x_{1} + \sum_{i=1}^{n}[(x_{i} -1)^{2}(1 + \sin^{2}(3 \pi x_{i+1}))])\}$ \\ $+ (x_{n}-1)^{2}[1 + \sin^{2}(2 \pi x_{n})]+ \sum_{i =1}^{n}u(x_{i}, 5,100.4)$ \\ $u(x_{i},a,k,m) = \begin{cases} {k(x_{i}-a)^{m}} & {x_{i} >a} \\ {0} & {-a < x_{i} < a} \\ {k(-x_{i}-a)^{m}} \end{cases} $} & $[-50. 50]^{n}$ \\ \hline
    Pinter & \makecell{$g_{8}(x) = \sum_{i=1}^{n}ix_{i}^{2} + \sum_{i=1}^{n}20 i \sin^{2}(x_{i-1}\sin(x_{i}) + \sin(x_{i+1}))$ \\ $+ \sum_{i=1}^{n}i\log_{10}(1 + i(x_{i-1}^{2}-2x{i}+ 3x_{i+1}- \cos(x_{i})+1)^{2}) , x_{0}= x_{n}, x_{n+1}= x_{i}$ } & $[-10,10]^{n}$ \\ \hline
    Quintic & $g_{9}(x) = \sum_{i=1}^{n} \vert x_{i}^{5} - 3x_{i}^{4} + 4x_{i}^{3}+ 2x_{i}^{2} - 10x_{i} -4  \vert $ & $[-10,10]^{n}$ \\ \hline
    Salomon & $g_{10}(x) = 1 - cos(2\pi \sqrt{\sum_{i=1}^{n}}) + 0.1\sqrt{\sum_{i=1}^{n}}$ & $[-100,100]^{n}$ \\ \hline
    sargan & $g_{11}(x) = \sum_{i=1}^{n} n(x_{i}^{2}+ 0.4\sum_{i \ne j}x_{i}x_{j})$ & $[-100,100]^{n}$ \\ \hline
    Sinusoidal & $g_{12}(x) = - [2.5\Pi_{i=1}^{n}\sin(x_{i}-30) + \Pi_{i=1}^{n}\sin(5(x_{i}-30))]$ & $[0,180]^{n}$ \\ \hline
    Styblinski - Tang & $g_{13}(x) = (-\frac{1}{2}n)\sum_{i=1}^{n}(x_{i}^{4} - 16x_{i}^{2}+ 5x_{i})$ & $[-5,5]^{n}$ \\ \hline \hline
    
    \end{tabular}
    \caption{MULTIMODAL TEST FUNCTIONS USED IN THE EXPERIMENTAL STUDY}
\end{table}

\begin{table}[H]
    \centering
    \begin{tabular}{lll}
    \hline \hline
    Name & Function & Interval \\ \hline 
    Alpine 1 & $h_{1}(x) = \sum_{i=1}^{n} \vert  x_{i}\sin(x_{i}+ 0.1x_{i}) \vert $ & $[-10,10]^{n}$  \\ \hline
    Bent Cigar & $h_{2}(x) = x_{1}^{2} + 10^{6}\sum_{i=2}^{n}x_{i}^{2}$ & $[-100,100]^{n}$ \\ \hline
    Bohachevsky & $h_{3}(x) = \sum_{i=1}^{n-1}[ x_{i}^{2} + 2x_{i+1}^{2} - 0.3\cos(3 \pi x_{i}) -0.4\cos(4 \pi x_{i+1}) + 0.7 ]$ & $[-50,50]^{n}$ \\ \hline
    Chung Reynold & $g_{4}(x) = (\sum_{i=1}^{n}x_{i}^{2})^{2}$ & $[-100,100]^{n}$ \\ \hline
    Csendes & $h_{5}(x) = \sum_{i=1}^{n} x_{i}^{6}(2+ \sin(\frac{1}{x_{i}}))$ & $[-1,1]^{n}$ \\ \hline
    Ellipsoid & $h_{6}(x) = \sum_{i=1}^{n}(100^{\frac{i-1}{n-1}}x_{i})^{2}$ &  $[-100,100]^{n}$ \\ \hline
    Rastrigin & $h_{7}(x)= \sum_{i=1}^{n}[x_{i}^{2}-10\cos(2\pi x_{i}) + 10]$ & $[-5.12, 5.12]^{n}$ \\ \hline
    Sum of different powers & $h_{8}(x) = \sum_{i=1}^{n}\vert x_{i} \vert^{i+1}$ & $[-1,1]^{n}$ \\ \hline
    Trigonometric 1 & $h_{9} = \sum_{i=1}^{n}[ n- \sum_{j=1}^{n}\cos(x_{j}) + i(1- \cos(x_{i}) - \sin(x_{i})) ]^{2}$ & $[0, \pi]^{n}$\\ \hline
    Wavy & $h_{10}(x) = 1 - \{ (\frac{1}{n}) \sum_{i=1}^{n}[\cos(10 x_{i}) \exp{\frac{-x_{i}^{2}}{2}}] \}$ & $[-\pi, \pi]^{n}$ \\ \hline
    zakharov & $h_{11}(x) = \sum_{i=1}^{n} x_{i}^{2} + (0.5 \sum_{i=1}^{n}ix_{i})^{2} + (0.5\sum_{i=1}^{n}ix_{i})^{4} $ & $[-5,10]^{n}$ \\ \hline \hline    
    \end{tabular}
    \caption{BASIC FUNCTIONS WHOSE SHIFTED VERSIONS HAVE BEEN USED IN THE EXPERIMENTAL STUDY}
\end{table}

\section{References}

[1]	Alba, E., Dorronsoro, B., “The exploration/exploitation tradeoff in dynamic cellular genetic algorithms”, IEEE Transactions on Evolutionary Computation, vol. 9, pp. 126–142, 2005.

[2]	Bazaraa, M. S., Sherali, H. D., Shetty, C. M., “Nonlinear programming: Theory and Algorithms”, John Wiley and Sons, New York, NY, 2006. 

[3]	Bhargava, V., Fateen, S. E. K., “A. Bonilla-Petriciolet, Cuckoo search: a new nature-inspired optimization method for phase equilibrium calculations”, Fluid Phase Equilibria, vol. 337, pp. 191-200, 2013. 

[4]	Cai, Z., Wang, Y., “A multiobjective optimization based evolutionary algorithm for constrained optimization”, IEEE Transactions on Evolutionary Computation, vol. 10, pp. 658-675, 2006. 

[5]	Changdar, C., Mahapatra, G. S., Kumar Pal, R., “An efficient genetic algorithm for multi-objective solid travelling salesman problem under fuzziness”, Swarm and Evolutionary Computation, vol. 15, pp. 27-37, 2014. 

[6]	Chen, Z., Wang, C., “Modelling RFID signal distribution based on neural network combined with continuous ant colony optimization”, Neurocomputing, vol. 123, pp. 354-361, 2014. 

[7]	Childs, J. Stoeber, H., J., “Self-oriented, other-oriented and socially prescribed perfectionism in employees: Relationships with burnout and engagement”, Journal of Workplace Behavioral Health, vol. 25, pp. 269-281, 2010. 

[8]	Cuevas, E., “A swarm optimization algorithm inspired in the behavior of the social-spider”, Expert Systems with Applications, vol. 40, no. 16, pp.  6374-6384, 2013. 

[9]	Dasgupta, D., Zbigniew, M., “Evolutionary algorithms in engineering applications”, Springer Science \& Business Media, 2013. 

[10]	DorigoM., Gambardella, L., Middendorf, M., Stutzle, T., “Special section on ant colony optimization”, IEEE Transactions on Evolutionary Computation, vol. 6, no. 4, pp. 317-365, 2002. 

[11]	Du, W., Li, B., “Multi-strategy ensemble particle swarm optimization for dynamic optimization”, Information Sciences, vol. 178, pp.  3096-3109, 2008. 

[12]	Eita, M. A., Fahmy, M. M., “Group counseling optimization”, Applied Soft Computing, vol. 22, pp. 585–604, 2014. 

[13]	Flett, G. L., Hewitt, P. L., Blankstein, K. R., Solnik, M., Van Brunschot, M., “Perfectionism, social problem-solving ability, and psychological distress”, Journal of Rational-Emotive and Cognitive-Behavior Therapy, vol. 14, pp. 245-274, 1996. 

[14]	Formato, R. A., “Central force optimization: a new nature inspired computational framework for multidimensional search and optimization”, Studies in Computational Intelligence, vol. 129, pp. 221-238, 2008.  

[15]	Garai, G., Chaudhurii, B. B., “A novel hybrid genetic algorithm with Tabu search for optimizing multi-dimensional functions and point pattern recognition”, Information Sciences, vol. 221, pp. 28-48, 2013. 

[16]	Geem, Z. W., Kim, J. H., Loganathan, G. V., “A new heuristic optimization algorithm: harmony search”, Simulation, vol. 76, no. 2, pp. 60-68, 2001. 

[17]	Ghaemi, M., Feizi-Derakhshi, M. R., “Forest optimization algorithm”, Expert Systems with Applications, vol. 41, no. 15, pp. 6676-6687, 2014. 

[18]	Ghodousian, A., Raeisian Parvari, M., “A modified PSO algorithm for linear optimization problem subject to the generalized fuzzy relational inequalities with fuzzy constraints (FRI-FC)”, Information Sciences,  vol. 418-419, pp. 317-345, 2017. 

[19]	Ghodousian, A., Babalhavaeji, A., “An efﬁcient genetic algorithm for solving nonlinear optimization problems deﬁned with fuzzy relational equations and max-Lukasiewicz composition”, Applied Soft Computing, vol. 69, pp. 475–492, 2018. 

[20]	Ghodousian, A., Naeeimib, M., Babalhavaeji, A., “Nonlinear optimization problem subjected to fuzzy relational equations deﬁned by Dubois-Prade family of t-norms”, Computers \& Industrial Engineering, vol. 119, pp. 167–180, 2018. 

[21]	Ghorbani, N., Babaei, E., “Exchange market algorithm”, Applied Soft Computing, vol. 19, pp. 177–187, 2014. \

[22]	Hatamlou, A., “Black hole: A new heuristic optimization approach for data clustering”, Information Sciences, vol. 222, pp. 175-184, 2013. 

[23]	He, S., Wu, Q. H., Saunders, J., “Group search optimizer: an optimization algorithm inspired by animal searching behavior”, IEEE Transactions on Evolutionary Computation, vol. 13, pp. 973–990, 2009. 

[24]	Hewitt, P. L., Flett, G. L., “Perfectionism in the self and social contexts: Conceptualization, assessment, and association with psychopathology”, Journal of Personality and Social Psychology, vol. 3, pp. 456-470, 1991. 

[25]	Hewitt, P. L., Flett, G. L., Perfectionism and stress processes in psychopathology, In Flett, G. L., and Hewitt, P. L., (Eds.), Perfectionism: Theory, research, and treatment (pp. 255-284), Washington, DC: American Psychological Association, 2002. 

[26]	Hewitt, P. L., Flett, G. L., “Multidimensional Perfectionism Scale (MPS): Technical manual”, Toronto: Multi-Health Systems, 2004. 

[27]	Holland, J. H., “Adaptation in natural and artificial systems: An introductory analysis with applications to biology, Control, and Artificial Intelligence”, U Michigan Press, 1975.  

[28]	Jamil, M., Yang, X. S., “A literature survey of benchmark functions for global optimization problems”, International Journal of Mathematical Modelling and Numerical Optimization, vol. 4, no. 2, pp. 150-194, 2013.  

[29]	Kennedy, J., Eberhart, R.C., “Particle swarm optimization”, Proceedings of IEEE International Conference on Neural Networks, vol. 4, pp. 1942-1948, 1995. 

[30]	Kern, S., Müller, S. D., Hansen, N., Büche, D., ocenasek, J., Koumoutsakos, P., “Learning probability distributions in continuous evolutionary algorithms-A comparative review”, Natural Computing, vol. 3, no. 1, pp. 77-112,  2004.  

[31]	Kim, D. H., Abraham, A., Cho, J. H., “A hybrid genetic algorithm and bacteria foraging approach for global optimization”, Information Sciences, vol. 177, pp. 3918-3937, 2007. 

[32]	Kirkpatrick, S., Gelatto, C. D., Vecchi, M. P., “Optimization by simulated annealing”, Science, vol. 220, pp. 671-680, 1983. 

[33]	Malviya, R., Pratihar, D. K., “Tuning of neural networks using particle swarm optimization to model MIG welding process”, Swarm and Evolutionary Computation, vol. 1, no. 4, pp. 223-235, 2011. 

[34]	Mirjalili, S., Mirjalili, S. M., Lewis, A., “Grey wolf optimizer”, Advances in Engineering Software, vol. 69, pp. 46-61, 2014. 

[35]	Mirjalali, S., Lewis, A., “The whale optimization algorithm”, Advances in Engineering Software, vol. 95, pp. 51-67, 2016. 

[36]	Montes, E. M., Coello, C., “A simple multimembered evolution strategy to solve constrained optimization problems”, IEEE Transactions on Evolutionary Computation, vol. 9, pp. 1-17, 2005. 

[37]	Nanda, S. J., Panda, G., “A survey on nature inspired methaheuristic algorithms for partitional clustering”, Swarm and Evolutionary Computation, vol. 16, pp.1-18, 2014. 

[38]	Pham, D. T., Ghanbarzadeh, A., Koc, E., Otri, S., Rahim, S., Zaidi, M., “The Bees Algorithm-A novel tool for complex optimization problems”, Manufacturing Engineering Centre, Cardiff University, Cardiff CF24 3AA, UK, 2005.  

[39]	Rajabioun, R., “Cuckoo optimization algorithm”, Applied Soft Computing, vol. 11, pp. 5508-5518, 2011. 

[40]	Rao, R. V., Savsani, V. J., Vakharia, D. P., “Teaching–learning-based optimization: an optimization method for continuous non-linear large scale problems”, Information Sciences, vol. 183, pp. 1–15, 2012. 

[41]	Rashedi, E., Nezamabadi-pour, H., Saryazdi, S., “GSA: A gravitational search algorithm”, Information Sciences, vol. 179, pp. 2232-2248, 2009.

[42]	Rechenberg, I., “Evolutionsstrategien”, Springer Berlin Heidelberg, 1978 pp. 83–114. 

[43]	Reynolds, R.G., Peng, B., “Cultural Algorithms: computation modeling of how cultures learn to solve problems: an engineering example”, Cybernetics and Systems, vol. 36, pp. 753-771, 2005.  

[44]	Simon, D., “Biogeography-based optimization”, IEEE Transactions on Evolutionary Computations, vol. 12, no. 6, pp. 702-713, 2008. 

[45]	Smith, M. M., Sherry, S. B., Chen, S., Flett, G. L., Hewitt, P. L., “Perfectionism and narcissism: A meta-analytic review”, Journal of Personality Research, vol. 64, pp. 90-101, 2016. 

[46]	Smith, M. M., Sherry, S. B., Rnic, K., Saklofske, D. H., Enns, M., Gralnick, T., “Are perfectionism dimensions vulnerability factors for depressive symptoms after controlling for neuroticism? A meta-analysis of 10 longitudinal studies”, European Journal of Personality, vol. 30, pp. 201-212, 2016.   

[47]	Smith, M. M., Speth, T. A., Sherry, S. B., Saklofske, D. H., Stewart, S. H., Glowacka, M., “Is socially prescribed perfectionism veridical? A new take on the stressfulness of perfectionism”, Personality and Individual Differences, vol. 110, pp. 115-118, 2017.  

[48]	Socha, K., Dorigo, M., “Ant colony optimization for continuous domain”, European Journal of operational Research, vol. 185, pp. 1155-1173, 2008. 

[49]	Stoeber, J., “How other-oriented perfectionism differs from self-oriented and socially prescribed perfectionism: Further findings”, Journal of Psychopathology and Behavioral Assessment, vol. 37, pp. 611-623, 2015. 

[50]	Stoeber, J., “How other-oriented perfectionism differs from self-oriented and socially prescribed perfectionism”, Journal of Psychopathology and Behavioral Assessment, vol. 36, pp. 329-338, 2014.

[51]	Stoeber, J., Corr, P. J., Smith, M. M., Saklofske, D. H., Perfectionism and personality. In J. Stoeber (Ed.), The psychology of perfectionism: Theory, research,  applications. London: Routledge, 2018, pp. 68-88. 

[52]	Stoeber, J., Hoyle, A., Last, F., “The consequences of perfectionism scale: Factorial structure and relationships with perfectionism, performance perfectionism, affect, and depressive symptoms”, Measurement and Evaluation in Counseling and Development, vol. 46, no. 3, pp. 178-191, 2013. 

[53]	Storn, R.M., Price, K., “Differential evolution-A simple and efficient heuristic for global optimization over continuous spaces”, Journal of Global Optimization, vol. 11, pp. 341-359, 1997. 

[54]	Storn, R., “System design by constraint adaptation and differential evolution”, IEEE Transactions on Evolutionary Computation, vol. 3, pp. 22-34, 1999. 

[55]	Stroud, P., “Kalman-extended genetic algorithm for search in nonstationary environments with noisy fitness evaluations”, IEEE Transactions on Evolutionary Computation, vol. 5, pp. 66-77, 2001. 

[56]	Suresh, K., Kumarappan, N., “Hybrid improved binary particle swarm optimization approach for generation maintenance scheduling problem”, Swarm and Evolutionary Computation, vol. 9, pp. 69-89, 2013.  

[57]	Tarasewich, P., McMullen, P. R., “Swarm intelligence: power in numbers”, Communication of ACM, vol. 45, pp. 62-67, 2002.  

[58]	Wang, B., Jin, X., Cheng, B., “Lion pride optimizer: an optimization algorithm inspired by lion pride behavior”, Science China Information Sciences, vol. 55, no. 10, pp. 2369-2389, 2012. 

[59]	Wolpert, D. H., Macready, W. G., “No free lunch theorems for optimization”, IEEE transactions on evolutionary computation, vol. 1, no. 1, pp. 67-82, 1997.

[60]	Yang, X. S., “A new metaheuristic bat-inspired algorithm”, Proceedings of the workshop on nature inspired cooperative strategies for optimization (NICSO 2010), Springer, p. 65-74, 2010. 

[61]	Yang, X. S., “Firefly algorithm Stochastic Test Functions and Design Optimization”, Int. J. Bio-Inspired Computation, vol. 2, no. 2, pp. 78-84, 2010.

[62]	Yang, H., Stoeber, J., “The Physical Appearance Perfectionism Scale: Development and Preliminary Validation”, Journal of Psychopathology and Behavioral Assessment, vol. 34, no. 1, pp. 69–83, 2012. 

[63]	Yao, X., Liu, Y., Lin, G., “Evolutionary programming made faster”, IEEE Transactions on Evolutionary Computation, vol. 3, pp. 82-102, 1999. 

[64]	Yeh, W. C., “Novel swarm optimization for mining classification rules on thyroid gland data”, Information Sciences, vol. 197, pp. 65-76, 2012.

[65]	Zheng, Y. J., “Water wave optimization: a new nature-inspired metaheuristic”, Computers and Operations Research, vol. 55, pp. 1-11, 2015.

\end{document}